\documentclass[12pt]{amsart}
\usepackage{amssymb,amscd,amsthm,amsmath,color}
\usepackage{fullpage}
\usepackage{epsfig}
\usepackage{graphicx}
\usepackage{xypic}

\numberwithin{equation}{section}

\newcommand{\PP}{\mathbb{P}}
\newcommand{\OO}{\mathcal{O}}

\newcommand{\FF}{\mathcal{F}}

\newcommand{\gG}{\mathbb{G}}

\newcommand{\cC}{\mathcal{C}}

\newcommand{\Tor}{\textsf{Tor}}
\newcommand{\rk}{\operatorname{rk}}
\newcommand{\degbd}{\operatorname{degbd}}

\newcommand{\Eff}{\operatorname{Eff}}
\newcommand{\Nef}{\operatorname{Nef}}
\newcommand{\Supp}{\operatorname{Supp}}

\newcommand{\Rat}{\operatorname{Rat}}

\newcommand{\Gr}{\operatorname{Gr}}
\newcommand{\SP}{\operatorname{SP}}
\newcommand{\ESP}{\operatorname{ESP}}
\newcommand{\Vol}{\operatorname{Vol}}

\newcommand{\Hilb}{\operatorname{Hilb}}
\newcommand{\Br}{\mathrm{Br}}

\newtheorem{theorem}{Theorem}[section]
\newtheorem{lemma}[theorem]{Lemma}
\newtheorem{proposition}[theorem]{Proposition}
\newtheorem{corollary}[theorem]{Corollary}

\theoremstyle{definition}
\newtheorem{definition}[theorem]{Definition}
\newtheorem{conjecture}[theorem]{Conjecture}
\newtheorem{remark}[theorem]{Remark}
\newtheorem{observation}[theorem]{Observation}

\newtheorem{example}[theorem]{Example}

\input xy
\xyoption{all}

\begin{document}

\title{Restricted tangent bundles for general free rational curves}

\author{Brian Lehmann}
\address{Department of Mathematics \\
Boston College  \\
Chestnut Hill, MA \, \, 02467}
\email{lehmannb@bc.edu}

\author{Eric Riedl}
\address{Department of Mathematics \\
University of Notre Dame  \\
255 Hurley Hall \\
Notre Dame, IN 46556}
\email{eriedl@nd.edu}

\begin{abstract}
Suppose that $X$ is a smooth projective variety and that $C$ is a general member of a family of free rational curves on $X$.  We prove several statements showing that the Harder-Narasimhan filtration of $T_{X}|_{C}$ is approximately the same as the restriction of the Harder-Narasimhan filtration of $T_{X}$ with respect to the class of $C$.  When $X$ is a Fano variety, we prove that the set of all restricted tangent bundles for general free rational curves is controlled by a finite set of data.  We then apply our work to analyze Peyre's ``freeness'' formulation of Manin's Conjecture in the setting of rational curves.  
\end{abstract}

\maketitle

\section{Introduction}

Let $X$ be a smooth projective variety over $\mathbb{C}$.  For any rational curve $f: \mathbb{P}^{1} \to X$ the restricted tangent bundle $f^{*}T_{X}$ decomposes into a sum of line bundles and this decomposition controls the deformation theory of $f$.  The possible direct summand decompositions of $f^{*}T_{X}$ have been extensively studied for special varieties including projective space \cite{Ascenzi88}, \cite{Ramella}, \cite{Ran01}, \cite{GHI13}, \cite{BGI16}, homogeneous varieties \cite{Mandal}, and hypersurfaces \cite{BrowningSawin}.  They have also been studied by Hwang and Mok in the context of minimal rational curves; see \cite{HwangMok, Hwang98} for more details.

Our goal is to study the behavior of the restricted tangent bundle for rational curves on arbitrary Fano varieties.  We will concentrate on free rational curves which are general in moduli -- the restricted tangent bundle of any other free curve will be a specialization.  Our main results show that for general curves the Harder-Narasimhan filtration of the restricted tangent bundle approximately corresponds to the Harder-Narasimhan filtration of the tangent bundle with respect to the curve class;  for highly multiple curve classes in dimension at most 5, the correspondence is exact. Furthermore, we show that for Fano varieties these filtrations are ``controlled'' by a finite set of data.  

\bigskip

Given a smooth projective variety $X$, we let $\Rat(X)$ denote the closure in $\overline{M}_{0,0}(X)$ of the sublocus parametrizing stable maps with irreducible domains onto free rational curves.  The following definition identifies our main object of interest.

\begin{definition}
\label{def-introExpectedSlopePanel}
Let $X$ be a smooth projective variety of dimension $n$ and let $f: \mathbb{P}^{1} \to X$ be a morphism with image curve $C = f_{*}\mathbb{P}^{1}$.  We will write $f^{*}T_{X} \cong \oplus_{i=1}^{n} \mathcal{O}(a_{i})$ where the $a_{i}$ are non-increasing.  We also let $\mu_{C}(T_{X}) = \frac{-K_{X} \cdot C}{n}$ denote the slope of $T_{X}$ with respect to $C$.
The \emph{slope panel} of the rational curve $C$ is the $n$-tuple
\begin{equation*}
\SP(C) := \left( \frac{a_{1}}{\mu_{C}(T_{X})}, \frac{a_{2}}{\mu_{C}(T_{X})}, \ldots \frac{a_{n}}{\mu_{C}(T_{X})}  \right).
\end{equation*}
\end{definition}

\vspace{.05in}

The primary factor in determining the slope panel of $C$ is the Harder-Narasimhan filtration of $T_{X}$ with respect to the class of $C$.  Loosely speaking, one expects the Harder-Narasimhan filtration of $T_{X}|_{C}$ to be ``close'' to the restriction of the Harder-Narasimhan filtration of $T_{X}$ with respect to the class of $C$.  In other words, the restriction of the Harder-Narasimhan filtration yields the ``expected'' value of the slope panel.

\begin{definition}
Let $X$ be a smooth projective variety of dimension $n$ and let $M$ be an irreducible component of $\Rat(X)$.  Let $\psi: Y \to X$ be a resolution of the Stein factorization of the evaluation map for the normalization of the irreducible component of $\overline{M}_{0,1}(X)$ lying over $M$.  Let $C'$ denote the image on $Y$ of a general curve parametrized by $M$ and denote the Harder-Narasimhan filtration of $T_{Y}$ with respect to the class of $C'$ by
\begin{equation*}
0 = \mathcal{F}_{0} \subset \mathcal{F}_{1} \subset \ldots \subset \mathcal{F}_{r} = T_{Y}.
\end{equation*}
For each $i \in \{1,\ldots,n\}$, we let $\sigma(i)$ denote the smallest integer such that $\rk(\mathcal{F}_{\sigma(i)}) \geq i$.  Set  $b_{i} = \mu_{C'}(\FF_{\sigma(i)}/\FF_{\sigma(i)-1})$.  We define the \emph{expected slope panel} of $C$ to be the $n$-tuple
\begin{equation*}
\ESP(C) := \left( \frac{b_{1}}{\mu_{C}(T_{X})}, \frac{b_{2}}{\mu_{C}(T_{X})}, \ldots \frac{b_{n}}{\mu_{C}(T_{X})}  \right).
\end{equation*}
(In other words, the slope panel has $\rk(\FF_{i}/\FF_{i-1})$ entries equal to $\mu_{C'}(\FF_{i}/\FF_{i-1})/\mu_{C}(T_{X})$.)
\end{definition}

\begin{remark}
Note that the expected slope panel for $C$ will be determined on a generically finite cover of $X$ when the evaluation map for the universal family over $M$ fails to have connected fibers.
\end{remark}

Our first theorem shows that the slope panel is usually very close to the expected slope panel.  Fix a positive integer $k$ and an irreducible component $M$ of $\Rat(X)$.  Choose an irreducible component of the sublocus of $M$ parametrizing stable maps through a general rational point $p$ in $X$, take $k$ general rational curves in this locus, and glue them at $p$ to obtain a stable map.  It turns out that this prescription uniquely determines an irreducible component of $\Rat(X)$ (see Definition \ref{defi:rgluedcomponents}).  As we vary the components $M$, we let $\Rat^{(k)}(X)$ denote the subset of $\Rat(X)$ obtained via this construction.

\begin{theorem} \label{theo:maintheorem1}
Let $X$ be a smooth projective uniruled variety of dimension $n$. 
\begin{enumerate}
\item Fix $\epsilon > 0$.  Then for any irreducible component $M$ of $\Rat(X)$ that parametrizes curves of anticanonical degree $\geq n^{2}/2\epsilon$ a general curve $C$ parametrized by $M$ will satisfy
\begin{equation*}
|\SP(C) - \ESP(C)|_{sup} < \epsilon.
\end{equation*}
\item Suppose $n \leq 5$.  Then there is some constant $k=k(n)$ (independent of $X$) such that a general curve $C$ parametrized by an irreducible component $M$ of $\Rat^{(k)}(X)$ satisfies $\SP(C) = \ESP(C)$.
\end{enumerate}
\end{theorem}

The first statement is a direct consequence of the Grauert-Mulich theorem recently proved in \cite[Proposition 3.1]{PatelRiedlTseng}.  The second statement relies on a ``Mehta-Ramanthan'' type result which we discuss in Section \ref{sect:intromr} below.

\begin{remark}
When $T_{X}$ is semistable we might expect that a ``sufficiently positive'' free rational curve will have restricted tangent bundle which is as close to semistable as possible in the sense that its direct summands differ by at most $1$.  However, \cite{Mandal} shows that this property fails already for the Grassmannian $\mathbb{G}(1,3)$; see Example \ref{ex-Grassmannian}.  Thus we can only hope for such a statement for a ``sufficiently divisible'' curve as in Theorem \ref{theo:maintheorem1}.(2).
\end{remark}

It is natural to ask how the expected slope panel changes as we vary the numerical class of the family of rational curves.  When $X$ is Fano, we can prove a strong finiteness statement that combines the results of \cite{Neumann10} with techniques from \cite{LST18}.

\begin{theorem} \label{theo:maintheorem2}
Suppose $X$ is a smooth Fano variety of dimension $n$.  There is a finite set of generically finite morphisms $\{ \phi_{i}: Y_{i} \to X\}$ from smooth projective varieties $Y_{i}$ and for every $i$ a finite set of homomorphisms $\{\ell_{ij}: N_{1}(Y_{i})_{\mathbb{Z}} \to \mathbb{Q}^{n}\}$ such that the following holds.

For every irreducible component $M$ of $\Rat(X)$ there is an index $i$ and an irreducible component $N$ of $\Rat(Y_{i})$ such that composition with $\phi_{i}$ induces a dominant rational map $\phi_{i*}: N \dashrightarrow M$.  Moreover, there is an index $j$ such that $\ESP(C) = \frac{\dim(X)}{-K_{X} \cdot C} \ell_{ij}(C')$ where $C,C'$ denote general curves parametrized by $M,N$ respectively.
\end{theorem}

\subsection{Peyre's formulation of Manin's Conjecture}
One of our motivations for this study is \cite{Peyre17} in which Peyre formulated a new version of Manin's Conjecture based on the notion of the ``freeness'' of a rational point.  Peyre suggests that the exceptional set in Manin's Conjecture over a number field should include all points with small ``freeness''.

\cite{LT19} formulates a version of Manin's Conjecture for rational curves known as Geometric Manin's Conjecture.  This conjecture predicts the behavior of the discrete invariants of irreducible components of $\Rat(X)$ -- the dimension and the number of components -- as the degree increases.  Just as in Manin's Conjecture for rational points, one must remove an ``exceptional set'' of components to obtain the expected asymptotic behavior.  We will use a construction of the exceptional set that is inspired by (but not identical to) the definitions of \cite{LT19}, \cite{LST18}.  Our goal is to contrast this construction with Peyre's formulation.

By analogy with \cite{Peyre17}, we make the following definition:

\begin{definition}[\cite{Peyre17} D\'efinition 4.11]
The \emph{minimal slope ratio} of a free curve $C$ is the smallest entry $a_{n}/\mu_{C}(T_{X})$ in the slope panel of $C$.
\end{definition}

It is then interesting to ask how the exceptional set in Geometric Manin's Conjecture interacts with the set of irreducible components of $\Rat(X)$ which have small minimal slope ratio.    The example of \cite{Sawin20} shows that Peyre's formulation may not ``remove enough'' -- there can be components with large minimal slope ratio which dominate the asymptotic count.  On the other hand, our work suggests that Peyre's formulation will never ``remove too much.''

More precisely, given a Fano variety $X$ Theorem \ref{theo:maintheorem1} constrains the set of lattice points in $\Nef_{1}(X)$ which can be represented by components of $\Rat(X)$ with small minimal slope ratio.  Assuming that irreducible components of $\Rat(X)$ are ``equidistributed'' in $\Nef_{1}(X)$ (see Conjecture \ref{conj:gmc}) this implies that the set of components removed by \cite{Peyre17} but not by \cite{LST18} will only make a negligible contribution to the asymptotic counting function.

In the following statement, $N(X,-K_{X},q,d)$ denotes the counting function in Geometric Manin's Conjecture (reflecting the dimension and number of components of $\Rat(X)$) and $N^{\ell > \epsilon}(X,-K_{X},q,d)$ denotes the modified counting function of \cite{Peyre17}.

\begin{theorem}
Assume Geometric Manin's Conjecture (Conjecture \ref{conj:gmc}).  Let $X$ be a smooth projective Fano variety and let $\epsilon: [1,\infty) \to (0,1)$ be a continuous decreasing function such that $\lim_{d \to \infty} \epsilon(d) = 0$.  For any $\delta > 0$, there is some $d_{0}$ sufficiently large such that
\begin{equation*}
N^{\ell > \epsilon}(X,-K_{X},q,d) > (1-\delta) N(X,-K_{X},q,d)
\end{equation*}
for every $d \geq d_{0}$.
\end{theorem}

We will define the various terms and state a precise version of Geometric Manin's Conjecture in Section \ref{sect:gmc}.  

\subsection{Deformations of vector bundles on rational curves} \label{sect:intromr}

Our main tool for relating the behavior of Harder-Narasimhan filtrations to restricted tangent bundles is the version of the Grauert-Mulich theorem recently proved in \cite[Proposition 3.1]{PatelRiedlTseng}.  For vector bundles of low rank, we will prove an even stronger statement:

\begin{theorem} \label{theo:mrforfreecurves}
Let $X$ be a smooth projective variety and let $\mathcal{E}$ be a torsion-free coherent sheaf on $X$ of rank $\leq 5$.  There is some constant $r$ (independent of $X$ and $\mathcal{E}$) such that the following holds.

Let $M$ be a component of $\Rat(X)$ such that the evaluation map for a resolution of the universal family over $M$ has connected fibers and let $C$ be a general curve parametrized by $M$.  If $\mathcal{E}$ is semistable with respect to the class of $C$, then there is a free rational curve $C'$ obtained by gluing $r$ members of $M$ and smoothing such that $\mathcal{E}|_{C'}$ is semistable.
\end{theorem}

\begin{remark}
This result is reminiscent of the Mehta-Ramanathan theorem which addresses complete intersection classes of the form $H^{n-1}$ for an ample divisor $H$.  We do not know whether an analogous ``Mehta-Ramanathan result'' for free rational curves holds when the torsion free sheaf $\mathcal{E}$ has arbitrary rank.
\end{remark}

Our proof of Theorem \ref{theo:mrforfreecurves} has two inputs.  First, we perform a careful analysis of how vector bundles behave on a smooth family of rational curves specializing to a union of two curves meeting at a node. Theorem \ref{theo:degreebound}  bounds the possible direct summands of the restriction of a bundle to a general curve in terms of its restriction to the nodal curve.   Second, given a vector bundle $\mathcal{E}$ on a projective variety $X$ we study when the positive summands of the restriction of $\mathcal{E}$ to a family of rational curves can be combined to yield a saturated subsheaf.

\bigskip

\noindent
{\bf Acknowledgements:}
The authors thank Izzet Coskun for suggesting the Quot scheme and for several helpful comments about Grassmannians.  We also thank Anand Deopurkar for an example clarifying the importance of smoothness in Section \ref{sect:vectorbundles}.  We are grateful to Jason Starr for helpful comments on \cite{Kollarbook}.  Part of this work was done at an AIM SQuaRE workshop, and we thank AIM for the excellent working environment.  Brian Lehmann is supported by NSF grant 1600875. Eric Riedl is supported by NSF grant 1945144.

\section{Preliminaries}

Throughout we will work with separated finite type schemes over $\mathbb{C}$. All varieties are reduced and irreducible.  For a projective variety $X$ we will use $\Eff^{1}(X) \subset N^{1}(X)$ to denote the pseudo-effective cone of divisors and $\Nef_{1}(X) \subset N_{1}(X)$ to denote the nef cone of curves.

\subsection{Vector bundles on rational curves}

Let $\mathcal{E}$ denote a rank $r$ degree $d$ vector bundle on $\mathbb{P}^{1}$.  We can write
\begin{equation*}
\mathcal{E} \cong \oplus_{i=1}^{r} \mathcal{O}(a_{i})
\end{equation*}
where $a_{i} \geq a_{j}$ for $i < j$ and $\sum_{i=1}^{r} a_{i} = d$.  We say that $\mathcal{E}$ is $k$-balanced if $\sup_{i,j} |a_{i} - a_{j}| = k$.  Note that if $\mathcal{E}$ is $0$-balanced or $1$-balanced then $\mathcal{E}$ is ``as balanced as possible'' given its degree and rank.  

We will need the following result showing that vector bundles  on $\mathbb{P}^{1}$ can only become ``more balanced'' under generalization.

\begin{lemma}[\cite{3264} Theorem 14.7.(a)] \label{lemm:deformingonp1} 
Let $\mathcal{E}_{1}, \mathcal{E}_{2}$ be two vector bundles on $\mathbb{P}^{1}$ of equal rank $r$ and degree $d$, and write
\begin{equation*}
\mathcal{E}_{1} = \oplus_{i=1}^{r} \mathcal{O}(a_{i}) \qquad \qquad \mathcal{E}_{2} = \oplus_{j=1}^{r} \mathcal{O}(b_{j})
\end{equation*}
where the $\{a_{i}\}$ and $\{b_{j}\}$ form non-increasing sequences.  The vector bundle $\mathcal{E}_{1}$ specializes to $\mathcal{E}_{2}$ if and only if for every integer $k$ satisfying $1 \leq k \leq r$ we have
\begin{equation*}
\sum_{i=r-k+1}^{r} a_{i} \geq \sum_{j=r-k+1}^{r} b_{j}.
\end{equation*}
\end{lemma}

\begin{definition}
Suppose that $\mathcal{E}$ is a vector bundle on $\mathbb{P}^{1}$ of the form $\mathcal{E} \cong \oplus_{i=1}^{r} \mathcal{O}(a_{i})$ where the $a_{i}$ are non-increasing.  Set $\mu = \deg(\mathcal{E})/r$.  The \emph{slope panel} of $\mathcal{E}$ is the $r$-tuple of rational numbers
\begin{equation*}
\left( \frac{a_{1}}{\mu}, \ldots, \frac{a_{r}}{\mu} \right).
\end{equation*}
The \emph{minimal slope ratio} is the smallest entry in the slope panel.
\end{definition}

\begin{definition}
We will say that a vector bundle $\mathcal{E}$ on $\mathbb{P}^{1}$ is \emph{sequential} if when we write $\mathcal{E} \cong \oplus_{i=1}^{r} \mathcal{O}(a_{i})$ as above we have $a_{i} \geq a_{i-1} - 1$ for $1 < i \leq r$.
\end{definition}

\subsection{Saturated subsheaves}

Let $X$ be a smooth projective variety equipped with a locally free sheaf $\mathcal{E}$.  In this section we briefly recall some facts about saturated subsheaves of $\mathcal{E}$.

\begin{lemma} \label{lemm:extendingopensetsat}
Let $X$ be a smooth projective variety equipped with a locally free sheaf $\mathcal{E}$ and let $U \subset X$ be an open subset.  Any saturated subsheaf $\mathcal{F}_{U} \subset \mathcal{E}|_{U}$ extends in a unique way to a saturated subsheaf $\mathcal{F} \subset \mathcal{E}$.
\end{lemma}

\begin{proof}
Consider the subsheaf $\mathcal{F}_{U} \subset \mathcal{E}|_{U}$.  By combining \cite[Lemma 28.22.2]{stacks-project} with the fact that a finite type subsheaf of a coherent sheaf is coherent, we see that there is a coherent subsheaf $\mathcal{F}' \subset \mathcal{E}$ such that $\mathcal{F}'|_{U} = \mathcal{F}_{U}$.  Let $\mathcal{F}$ denote the saturation of $\mathcal{F}'$ in $\mathcal{E}$.  Since $\mathcal{F}_{U}$ is saturated we have $\mathcal{F}|_{U} = \mathcal{F}_{U}$.

We next argue that there is a unique such extension.  Suppose given two saturated subsheaves $\mathcal{F}_{1},\mathcal{F}_{2}$ of $\mathcal{E}$ extending $\mathcal{F}_{U}$.  Consider the composition of the inclusion $\mathcal{F}_{1} \to \mathcal{E}$ with the quotient $\mathcal{E} \to \mathcal{E}/\mathcal{F}_{2}$.  The composition vanishes along $U$, and since $\mathcal{E}/\mathcal{F}_{2}$ is torsion-free we conclude that the map vanishes everywhere.  As the reverse is also true by symmetry we see that $\mathcal{F}_{1} = \mathcal{F}_{2}$.  
\end{proof}

\begin{lemma} \label{lemm:distandgrass}
Let $X$ be a smooth projective variety equipped with a locally free sheaf $\mathcal{E}$ of rank $n$.  There is a bijection between saturated subsheaves $\mathcal{F} \subset \mathcal{E}$ of rank $d$ and rational maps $X \dashrightarrow \Gr(n-d,\mathcal{E})$ to the relative Grassmannian of $(n-d)$-dimensional quotients.
\end{lemma}

\begin{proof}
First suppose we have a saturated subsheaf $\mathcal{F} \subset \mathcal{E}$.  Let $U \subset X$ denote an open subset such that $\mathcal{E}/\mathcal{F}|_{U}$ is locally free.  Then we obtain a morphism $U \to \Gr(n-d,\mathcal{E})$ using the universal property of the Grassmannian.

Conversely, suppose given a rational map $\phi: X \dashrightarrow \Gr(n-d,\mathcal{E})$ defined on an open set $U$.  By taking the $\phi$-pullback of the universal subsheaf on $\Gr(n-d,\mathcal{E})$, we obtain a saturated subsheaf $\mathcal{F}_{U} \subset \mathcal{E}|_{U}$.  We can then extend this to a saturated subsheaf of $\mathcal{E}$ on all of $X$ using Lemma \ref{lemm:extendingopensetsat}.

By the uniqueness in Lemma \ref{lemm:extendingopensetsat} these two constructions are inverse to each other.
\end{proof}

We next use Lemma \ref{lemm:distandgrass} to construct saturated subsheaves from dominant families of curves.

\begin{lemma}
\label{lemm:positiveSubsheaf}
Let $\mathcal{E}$ be a torsion-free sheaf on a smooth projective variety $X$.  Let $g: \cC \to X$ be the evaluation map for a family of free rational curves parametrized by an open subset of a component of $\Rat(X)$ and let $\FF_0$ be a torsion-free subsheaf of $g^*\mathcal{E}$. Then there is a saturated subsheaf $\FF$ of $\mathcal{E}$ whose fiber at a general point $p$ is the span of $(\FF_0)_{(f,p)}$ as $f: C \to X$ ranges over all curves in $\cC$ passing through $p$.

Furthermore, for a general curve $C$ parametrized by $\cC$ we have an injection $\FF_{0}|_{C} \to g^{*}\FF|_{C}$.
\end{lemma}

\begin{proof}
Recall that any torsion-free sheaf is locally free away from a codimension $2$ subset.  By \cite[II.3.7 Proposition]{Kollarbook} a general member $C$ of a family of free rational curves will not intersect a fixed codimension $2$ subset.  Thus, if we choose a general member $C$ of our family of curves we can assume that $C$ is contained in the locally free locus of any finite set of torsion-free sheaves.

Let $V \subset \cC$ be the open subset for which $g^{*}\mathcal{E}/\FF_{0}$ is locally free.  We obtain a morphism $V \to \Gr(\rk(g^{*}\mathcal{E}/\FF_{0}),g^{*}\mathcal{E})$.  Since $g$ is dominant, there is an open subset $U$ of $X$ contained in $f(V)$.  We can define a map $U \to \Gr(d,\mathcal{E})$ sending a point $p$ in $U$ to the span of the planes parametrized by $g^{-1}(p)$.  By Lemma \ref{lemm:distandgrass} this gives a saturated subsheaf $\mathcal{F}$ of $\mathcal{E}$ which is equal to the span of $(\FF_{0})_{(f,p)}$ as we vary $f: C \to X$ over general curves in $\cC$ through $p$.  By taking a closure, we see that $\mathcal{F}$ satisfies the first statement.

We still must show the final statement.  By construction there is an open subset $U \subset X$ such that for a general curve $C$ in our family, the map $\FF_{0}|_{C \cap f^{-1}U} \to g^{*}\mathcal{E}|_{C \cap f^{-1}U}$ factors through $g^{*}\mathcal{F}|_{C \cap f^{-1}U}$.  Since $\mathcal{F}$ is a saturated subsheaf of $\mathcal{E}$ and since $C$ is general, it is still true that $g^{*}\mathcal{F}$ is a saturated subsheaf of $g^{*}\mathcal{E}$.  Thus the factoring of the map $\FF_{0}|_{C} \to g^{*}\mathcal{E}|_{C}$ through $g^{*}\mathcal{F}|_{C}$ over $f^{-1}(U)$ extends to a factoring over all of $C$.
\end{proof}

\subsection{Stability of the tangent bundle for Fano varieties}

Our main results show that the restricted tangent bundle for rational curves depends upon the stability properties of the tangent bundle.  When $X$ is a Fano variety, the stability of the tangent bundle has been extensively studied, particularly in regards to its relationship with K-stability.  In this section, we briefly review some examples which help illustrate the possible range of behaviors for Fano varieties.

Many Fano varieties of Picard rank $1$ -- complete intersections in projective space (\cite{PW95}), many homogeneous varieties (\cite{Boralevi12}), all examples of dimension $\leq 5$ (\cite{Hwang98}) -- have stable tangent bundle.  Previously it was expected that a Fano variety of Picard rank $1$ will always have a stable tangent bundle.  However, \cite{Kanemitsu19} gave several examples of smooth horospherical varieties of Picard rank $1$ which have unstable tangent bundle.  In particular, for such varieties $T_{X}$ is not stable with respect to any nef curve class.

This same phenomenon can occur for Fano varieties of higher Picard rank. 

\begin{example}
Fix a projective space $\mathbb{P}^{n}$ and choose a non-increasing sequence of positive integers $\{ a_{i} \}_{i=0}^{m}$ whose sum $d$ satisfies $d \leq n$.  Let $X$ be the projective bundle $\mathbb{P}_{\mathbb{P}^{n}}(\oplus_{i} \mathcal{O}(a_{i}))$ equipped with the map $\pi: X \to \mathbb{P}^{n}$.  Then $X$ is a Fano variety.  Let $\xi$ denote a section of $\mathcal{O}_{X/\mathbb{P}^{n}}(1)$ and let $H$ denote the pullback of a hyperplane on $\mathbb{P}^{n}$.  Then
\begin{equation*}
\Eff^{1}(X) = \langle H, \xi - a_{0}H \rangle.
\end{equation*}
If we set $\alpha = \xi^m \cdot H^{n-1}$ and $\ell = H^{n} \cdot \xi^{m-1}$ we get the following intersection matrix:

\begin{table}[h!] \centering
\begin{tabular}{c|c|c}
& H & $\xi$ \\ \hline
$\alpha$ & 1 & d \\ \hline
$\ell$ & 0 & 1
\end{tabular}
\end{table}

Thus the nef cone of curves is generated by the class $\ell$ of a line in a fiber in $\pi$ and the class $\beta = \alpha + (a_{0}-d) \ell$.

Consider the restriction of the exact sequence
\begin{equation*}
0 \to T_{X/\mathbb{P}^{n}} \to T_{X} \to \pi^{*}T_{\mathbb{P}^{n}} \to 0
\end{equation*}
to a general free rational curve $f: C \to X$ of class $x \beta + y \ell$.  The restriction of the leftmost term has degree $x(ma_{0} + a_{0} - d) + y(m+1)$ and the restriction of the rightmost term has degree $x(n+1)$.  Choose parameters so that $a_{0} = d = n \geq 3$ and $m \geq 2$.  Then the sheaf on the left has slope at least $1$ larger than the sheaf on the right so that $T_X|_C$ will be unbalanced for any rational nef curve $C$.  (In fact, if we ignore the rational curves contained in the fibers of $\pi$ then the difference between the highest and lowest summands of $T_{X}|_{C}$ will be at least $n-2$ for every other rational nef curve $C$.) 
\end{example}

\section{ Families of rational curves}

Let $X$ be a smooth projective variety.  We will parametrize rational curves on $X$ using the coarse moduli space  $\overline{M}_{0,0}(X)$ of the Kontsevich space of stable maps.  Note that if we choose a stable map $f: Z \to X$ representing a point of this moduli space, the isomorphism class of the restricted tangent bundle $f^{*}T_{X}$ is independent of the choice of representative $f$.  Given a stable map $f: Z \to X$ with image curve $C$, we will sometimes abusively use $C$ in place of $f$ and write $\mathcal{E}|_{C}$ in place of $f^{*}\mathcal{E}$ for coherent sheaves $\mathcal{E}$ on $X$.

We will let $\Rat(X)$ denote the closure of the sublocus of $\overline{M}_{0,0}(X)$ parametrizing free curves $f: \mathbb{P}^{1} \to X$.  In other words, any irreducible component $M$ of $\Rat(X)$ will generically parametrize free curves.  For such a component $M$, we will call the corresponding component $M'$ of $\overline{M}_{0,1}(X)$ the one-pointed family over $M$.  Note that $M'$ comes equipped with a flat map $p: M' \to M$ and an evaluation map $ev: M' \to X$.  We let $\Rat_{conn}(X) \subset \Rat(X)$ denote the sublocus consisting of those components $M$ which have the property that if we precompose $ev: M' \to X$ with a normalization map $\nu: M'^{\nu} \to M'$ the composition $ev \circ \nu$ has connected fibers.

\subsection{Gluing free curves} \label{sect:gluingcurves}

We will frequently deal with stable maps $f: Z \to X$ where $Z$ is the union of two rational curves meeting transversally at a point and the restriction of $f$ to each component of $Z$ is free.

\begin{lemma} \label{lemm:smoothingforfreecurves}
Let $X$ be a smooth projective variety.  Let $Z$ denote the union of two rational curves meeting transversally at a point.  Suppose that $f: Z \to X$ is a stable map such that the restriction of $f$ to each component of $Z$ is free.  Then there is a smooth surface $S$ equipped with a flat projective morphism to a curve $\pi: S \to T$, a point $0 \in T$, and a morphism $g: S \to X$ such that
\begin{enumerate}
\item the fiber of $\pi$ over $0$ is isomorphic to $Z$,
\item the restriction of $g$ to the fiber over $0$ is $f$, and
\item the restriction of g to a general fiber of $\pi$ is a free curve on $X$.
\end{enumerate}
\end{lemma}

\begin{proof}
Since the restriction of $f$ to each irreducible component of $Z$ is free, we see that $f$ will be a smooth point of $\overline{M}_{0,0}(X)$. Thus there is a smooth complete curve $T$ mapping finitely to $\overline{M}_{0,0}(X)$ and meeting the locus of irreducible curves such that we can pull back the universal family to $T$. This gives rise to a surface $S$ mapping to $X$ via a map $g$ and mapping to $T$ via a map $\pi$, with $\pi^{-1}(0)$ isomorphic to $Z$. Because the components of $Z$ are free, $S$ will be smooth along $Z$, and we can resolve the other singularities of $S$ to ensure that $S$ is smooth. By construction a general fiber of $S \to T$ is smooth and irreducible and since $f^* T_X$ is globally generated the restriction of $g^*T_X$ to a general fiber will be as well.
\end{proof}

Suppose we have two free curves $C_{1}$, $C_{2}$ on $X$ which meet at a point $p$.  By choosing a branch of each curve at $p$ we can construct a stable map from a reducible curve $Z$ which has two rational components meeting transversally at a single point.  Since both curves are free, the corresponding map $f: Z \to X$ is a smooth point of the space of stable maps.  We will (somewhat imprecisely) refer to this process as ``gluing and smoothing'' $C_{1}$ and $C_{2}$.  (Note that we allow $C_{1} = C_{2}$.)

\begin{lemma}
\label{lem-glueNonACovers}
Let $X$ be a smooth projective variety.  Suppose that $M_{1}$ and $M_{2}$ are two irreducible components of $\Rat_{conn}(X)$.
Then there is only one irreducible component $M$ of $\overline{M}_{0,0}(X)$ obtained by selecting a general point $p \in X$ and gluing and smoothing two general free curves through $p$ from $M_{1}$ and $M_{2}$.  Furthermore, $M$ is also an irreducible component of $\Rat_{conn}(X)$.  
\end{lemma}

\begin{proof}
Let $C_{1}$ and $C_{2}$ denote free curves parametrized by $M_{1}$ and $M_{2}$ which meet at a general point of $X$.   Note that the stable map obtained by gluing $C_1$ and $C_2$ is a smooth point of $\overline{M}_{0,0}(X)$.  Since our assumption implies that the unique component of $M'_{1} \times_{X} M'_{2}$ which dominates both factors under the projection maps is irreducible, there is a unique component $M$ of $\overline{M}_{0,0}(X)$ that contains generic points of $M'_{1} \times_{X} M'_{2}$ and this component generically parametrizes rational curves.

We still must show that the evaluation map for the normalization of $M'$ has connected fibers.  Suppose we fix a general point of $X$.  The fiber over $x$ of the evaluation map for the one-pointed family over $M'_{1} \times_{X} M'_{2}$ is irreducible.  Note that the monodromy action acts transitively on the irreducible components of fibers of $M'$ over $x$.  Since this action preserves the irreducible components of $\overline{M}_{0,0}(X)$ containing the two components of our reducible curve, we see that every component of the fiber of $M'$ over $x$ contains a curve parametrized by $M'_{1} \times_{X} M'_{2}$ and we deduce that the fiber of $M'$ over $x$ is irreducible.
\end{proof}

\subsection{Grauert-Mulich}
Recently, \cite{PatelRiedlTseng} proved a version of the Grauert-Mulich theorem which holds for arbitrary projective varieties.

\begin{lemma}[\cite{PatelRiedlTseng} Proposition 3.1]
\label{lemm:citedGrauert-Mulich}
Let $X$ be a smooth projective variety and let $\mathcal{E}$ be a torsion free sheaf on $X$.  Suppose that $M$ is an irreducible component of $\Rat_{conn}(X)$.  Let $C$ be a general free curve parametrized by $M$.  If $\mathcal{E}|_C = \bigoplus_i \OO(a_i)$ with $a_1 \geq a_2 \geq \dots \geq a_r$, then either $a_i \leq a_{i+1}+1$ or there is an index $j$ and a subsheaf $\FF \subset \mathcal{E}$ such that $\FF|_C = \OO(a_1) \oplus \dots \OO(a_j)$.
\end{lemma}

We will usually use the following consequence:

\begin{corollary} \label{coro:semistablesequential}
Let $X$ be a smooth projective variety and let $\mathcal{E}$ be a torsion free sheaf on $X$.  Suppose that $M$ is an irreducible component of $\Rat_{conn}(X)$ and let $C$ be a general free curve parametrized by $M$.  If $\mathcal{E}$ is semistable with respect to $C$ then $\mathcal{E}|_{C}$ is sequential. 
\end{corollary}

\section{Vector bundles on smoothings of rational curves} \label{sect:vectorbundles}

Suppose we have a family of curves $\pi: \mathcal{C} \to B$ whose general fiber $C$ is a rational curve and whose special fiber $Z$ is a union of two rational curves meeting at a single node.  Given a vector bundle $\mathcal{E}$ on $\mathcal{C}$, we analyze how information about $\mathcal{E}|_{Z}$ can be used to constrain the possibilities for $\mathcal{E}|_{C}$.

It turns out that there is an essential difference between the case when $\mathcal{C}$ is smooth and the case when $\mathcal{C}$ is singular at the node of $Z$.  In this section, we will focus on the situation when $\mathcal{C}$ is smooth, which is the correct setting for our intended applications.  When $\mathcal{C}$ is allowed to be singular, there is much less that one can say: by \cite{Smith20},  in general the only restriction one can place on $\mathcal{E}|_{C}$ comes from upper semicontinuity of cohomology groups.

The main result of this section is Theorem \ref{theo:degreebound} which shows that when $\mathcal{E}|_{Z}$ is ``close'' to balanced, then $\mathcal{E}|_{C}$ must also be ``close'' to balanced.  The idea is to use projectivity of the Quot scheme to limit how negative a quotient of $\mathcal{E}|_{C}$ can be.

\subsection{Sheaves on nodal curves}
We start with a brief review of sheaf theory on a union of two $\mathbb{P}^{1}$s.  Throughout this section $Z$ will denote the union of two rational curves $Z_{1},Z_{2}$ meeting transversally at a single node $p$.  

Let $\mathcal{F}$ be a sheaf on $Z$.  We say that $\mathcal{F}$ is torsion free if every subsheaf $\mathcal{E} \subset \mathcal{F}$ has support of dimension $1$.  By \cite[VII.1 Lemme 3]{Seshadri82} any sheaf $\mathcal{F}$ has a canonical sequence
\begin{equation*}
0 \to \mathcal{F}_{tors} \to \mathcal{F} \to \mathcal{F}_{tf} \to 0
\end{equation*}
where $\mathcal{F}_{tors}$ is supported on a dimension $0$ subset and $\mathcal{F}_{tf}$ is a torsion free quotient of $\mathcal{F}$.

\cite[VIII.1 Proposition 3]{Seshadri82} proves a local structure theorem for torsion free modules near a nodal point of an arbitrary curve.  Since we are working with two rational curves, this local result extends to a global one by exactly the same argument.

\begin{lemma}[\cite{Seshadri82} VIII.1 Proposition 3] \label{lemm:torsionfreedecomposition}
Any torsion free sheaf $\mathcal{F}$ on $Z$ admits a unique expression
\begin{equation*}
\mathcal{F} \cong \mathcal{G} \oplus \mathcal{H}_{1} \oplus \mathcal{H}_{2}
\end{equation*}
where $\mathcal{G}$ is a vector bundle on $Z$ and $\mathcal{H}_{1}, \mathcal{H}_{2}$ are respectively vector bundles on $Z_{1}, Z_{2}$.
\end{lemma}

We will also need to know about the structure of locally free sheaves on $Z$.  First of all, an invertible sheaf is determined up to isomorphism by its degree on the two components $Z_{1}, Z_{2}$.  We will denote the invertible sheaf whose degree on $Z_{1}$ is $a$ and whose degree on $Z_{2}$ is $b$ by $\mathcal{O}_{Z}(a,b)$.  In general every locally free sheaf is a direct sum of the $\mathcal{O}_{Z}(a_{i},b_{i})$:

\begin{lemma}[\cite{Ran07} Proposition 5.1]
\label{lemm:nodalsplitting}
Every locally free sheaf on $Z$ splits into a direct sum of line bundles.
\end{lemma}

We will later need one additional fact about torsion free sheaves on $Z$.

\begin{lemma} \label{lemm:gnonvanishing}
Let $g: \mathcal{O}_{Z}(a,b) \to \mathcal{F}$ be a morphism of torsion-free sheaves on $Z$ such that for the node $p$ the restriction $g|_{p}$ is not the zero map.  Set $\mathcal{F}'$ to be the cokernel of $g$.  Write
\begin{equation*}
\mathcal{F} \cong \mathcal{G} \oplus \mathcal{H}_{1} \oplus \mathcal{H}_{2} \qquad \qquad (\mathcal{F}')_{tf}  \cong \mathcal{G}' \oplus \mathcal{H}'_{1} \oplus \mathcal{H}'_{2}
\end{equation*}
as in Lemma \ref{lemm:torsionfreedecomposition}.
If $g$ is injective then
\begin{equation*}
0 \leq \rk(\mathcal{G}) -  \rk(\mathcal{G}') \leq 1.
\end{equation*}
If $g$ is not injective, then
\begin{equation*}
\rk(\mathcal{G}) - \rk(\mathcal{G}') =  0.
\end{equation*}
\end{lemma}

\begin{proof}
Let $s, t_{1}, t_{2}$ denote respectively $\rk(\mathcal{G})$, $\rk(\mathcal{H}_{1})$, $\rk(\mathcal{H}_{2})$ and let $s', t'_{1}, t'_{2}$ denote respectively $\rk(\mathcal{G}')$, $\rk(\mathcal{H}_{1}')$, $\rk(\mathcal{H}_{2}')$.

We first show that $\rk(\mathcal{G}) - \rk(\mathcal{G}') \geq 0$.  Let $h: \mathcal{F} \to \mathcal{G}'$ denote the composition of $g$ with the projection map $\mathcal{F}' \to \mathcal{G}'$.  For any direct summand $\mathcal{O}_{Z_{1}}(a)$ of $\mathcal{F}$ the map $\mathcal{O}_{Z_{1}}(a) \to \mathcal{G}'$ vanishes at the node, and similarly for $Z_{2}$.  Since $h|_{p}$ is surjective, we see that $\rk(\mathcal{G}) \geq \rk(\mathcal{G}')$.

To prove upper bounds on $\rk(\mathcal{G}) - \rk(\mathcal{G}')$, first suppose that $g$ is injective.  By comparing ranks at the generic points of $Z_{1}$ and $Z_{2}$ we see that
\begin{equation*}
s + t_{1} = 1 + s' + t'_{1} \qquad \qquad s+t_{2} = 1 + s' + t'_{2}.
\end{equation*}
Using the fact that $g|_{p}$ is not zero and comparing ranks at $p$, we see see that
\begin{equation*}
s+ t_{1} + t_{2} = 1 + s' + t'_{1} + t'_{2} + \chi(\mathcal{F}'_{tors}|_{p}).
\end{equation*}
Combining these equations we see that 
\begin{equation*}
s - s' = 1 - \chi(\mathcal{F}'_{tors}|_{p}).
\end{equation*}
In particular this shows that $s-s' \leq 1$.

The argument when $g$ is not injective is exactly analogous.  Since $\mathcal{F}$ is torsion-free $g$ must factor through the quotient map from $\mathcal{O}_{Z}(a,b)$ to either $\mathcal{O}_{Z_{1}}(a)$ or $\mathcal{O}_{Z_{2}}(b)$.  Without loss of generality we may assume it is the former.  By comparing ranks at the generic points of $Z_{1}$ and $Z_{2}$ we see that
\begin{equation*}
s + t_{1} = 1 + s' + t'_{1} \qquad \qquad s+t_{2} = s' + t'_{2}.
\end{equation*}
Using the fact that $g|_{p}$ is not zero and comparing ranks at $p$, we see see that
\begin{equation*}
s+ t_{1} + t_{2} = 1 + s' + t'_{1} + t'_{2} + \chi(\mathcal{F}'_{tors}|_{p}).
\end{equation*}
Combining these equations, we see that
\begin{equation*}
s - s' = -\chi(\mathcal{F}'_{tors}|_{p}).
\end{equation*}
This proves that $s-s' \leq 0$.
\end{proof}

\subsection{Deformations}

We next turn to the problem of analyzing how vector bundles interact with smoothings of a nodal curve.  As before $Z$ will denote a union of two smooth rational curves $Z_{1}, Z_{2}$ meeting transversally at a single node $p$.

\begin{lemma} \label{lemm:restrictionoftfsheaf}
Let $\pi: \mathcal{C} \to B$ be a smoothing of $Z$ with smooth total family and with $\dim(B) = 1$.  Let $\FF$ be a torsion free sheaf on $\cC$ whose restriction to a general fiber $C$ of $\pi$ is a vector bundle of rank $r$ and degree $d$.  Suppose we write
\begin{equation*}
(\mathcal{F}|_{Z})_{tf} = \mathcal{G} \oplus \mathcal{H}_{1} \oplus \mathcal{H}_{2}
\end{equation*}
where $\mathcal{G}$ is locally free on $Z$ of rank $s$, $\mathcal{H}_{1}$ is locally free of rank $t_{1}$ on $Z_{1}$, and $\mathcal{H}_{2}$ is locally free of rank $t_{2}$ on $Z_{2}$.
Then $t_{1} = t_{2}$ and
\begin{equation*}
\chi((\mathcal{F}|_{Z})_{tf}) \leq d + s.
\end{equation*}
\end{lemma}

\begin{proof}
Since $\FF$ is torsion free, by the discussion in \cite[Section 1]{Hartshorne80} there is an injection $\FF \to \FF^{**}$ with torsion cokernel.
Since $\FF$ and $\FF^{**}$ have the same first chern class, the cokernel is supported in codimension $2$.  Since $\FF^{**}$ is a reflexive sheaf on a smooth surface \cite[Corollary 1.4]{Hartshorne80} shows that $\FF^{**}$ is locally free.  Altogether, we have the following short exact sequence.
\[ 0 \to \FF \to \FF^{**} \to T \to 0 \]
where $T$ is supported at finitely many points.  After restriction to $Z$ we have the exact sequence
\[ 0 \to T_0 \to \FF|_Z \to \FF^{**}|_Z \to T|_Z \to 0 , \]
where $T_0 = \Tor^1(T, \OO_Z)$ is torsion.

Since $\FF$ is torsion free and $\pi$ is dominant, $\FF$ is also torsion free as an $\mathcal{O}_{B}$-module.  Thus $\FF$ is flat over $B$ by \cite[Lemma 15.22.10]{stacks-project}.  By comparing Hilbert polynomials, we see that the rank of $\FF$ at the generic points of $Z_{1}$ and $Z_{2}$ is $r$.  By comparing ranks at the generic points of $Z_{1}$ and $Z_{2}$ we see that $s + t_{1} = r = s+t_{2}$.  Thus from now on we simply write $t$ for $t_{1} = t_{2}$.

Since $\FF$ is flat over $B$, we have 
\begin{align*}
d+r = \chi(\FF|_Z) & = \chi((\mathcal{F}|_{Z})_{tf}) + \chi((\mathcal{F}|_{Z})_{tors}).
\end{align*}
To prove the statement, it suffices to verify that $s + \chi((\mathcal{F}|_{Z})_{tors}) \geq r$.

By restricting $\FF|_Z \to \FF^{**}|_Z$ to $p$, we obtain
\begin{equation*}
\FF|_p \to \mathcal{O}_{p}^{\oplus r}.
\end{equation*}
We also have the surjection
\begin{equation*}
\FF|_{p} \to ((\FF|_{Z})_{tf})|_{p} \cong \mathcal{O}_{p}^{\oplus s+ 2t}.
\end{equation*}
Thus, we see that the kernel $K$ of $\FF|_{p} \to \FF^{**}|_{p}$ has length at least $s+2t - r = t$.  Since $T_{0}$ surjects onto the kernel of $\FF|_{Z} \to \FF^{**}|_{Z}$, we see that $K$ is a quotient of $T_{0}$.  Since $T_{0}$ is a torsion subsheaf of $\mathcal{F}|_{Z}$, we obtain
\begin{equation*}
s + \chi((\mathcal{F}|_{Z})_{tors}) \geq s + \chi(T_{0}) \geq s + t = r
\end{equation*}
as desired. 
\end{proof}

\begin{lemma} \label{lemm:surjectiononnodalcurve}
Let $\mathcal{E}$ be a vector bundle of rank $r$ on the nodal curve $Z = Z_1 \cup Z_2$ with $Z_i$ isomorphic to $\PP^1$.  Suppose there is a surjection $g: \mathcal{E} \to \mathcal{F}$ where $\mathcal{F}$ is a torsion free sheaf of the form
\begin{equation*}
\mathcal{F} = \mathcal{G} \oplus \mathcal{H}_{1} \oplus \mathcal{H}_{2}
\end{equation*}
as in Lemma \ref{lemm:torsionfreedecomposition}.  Then there are non-negative integers $j$, $k_{1}$, and $k_{2}$ and a direct summand $\mathcal{E}' \cong \oplus_{i=1}^{j+k_{1}+k_{2}} \mathcal{O}_{Z}(a_{i},b_{i})$ of $\mathcal{E}$ such that $g|_{\mathcal{E}'}$ is generically surjective,
\begin{equation*}
j+k_{1} = \rk(\mathcal{G}) + \rk(\mathcal{H}_{1}), \qquad \qquad j+k_{2} = \rk(\mathcal{G}) + \rk(\mathcal{H}_{2}),
\end{equation*}
and
\begin{equation*}
\rk(\mathcal{G}) + \sum_{i=1}^{j} (a_{i} + b_{i}) + \sum_{i=j+1}^{j+k_{1}} (a_{i}+1) + \sum_{i=j+k_{1}+1}^{j+k_{1}+k_{2}} (b_{i}+1) \leq \chi(\mathcal{F}).
\end{equation*}
\end{lemma}

\begin{proof}  
Let $s, t_{1}, t_{2}$ denote respectively $\rk(\mathcal{G})$, $\rk(\mathcal{H}_{1})$, $\rk(\mathcal{H}_{2})$ and let $m$ denote the sum $2s+t_{1}+t_{2}$.  We prove the statement by induction on $m$, where $m$ is allowed to decrease by $1$ or $2$ each time.

We start with the base case when $m=1$ or $m=2$.  There are several cases:
\begin{enumerate}
\item $s=0$, $t_{2} = 0$, and $0 \leq t_{1} \leq 2$.  Let $\mathcal{O}_{Z}(a,b)$ be a summand of $\mathcal{E}$ such that the map to $\mathcal{F}$ does not vanish.  Then the restriction of $g$ to $\mathcal{O}_{Z}(a,b)$ will factor through $\mathcal{O}_{Z_{1}}(a)$.  Let $\mathcal{F}'$ denote the cokernel of $\mathcal{O}(a,b) \to \mathcal{F}$.  Then
\begin{equation*}
\chi(\mathcal{F}) - \chi(\mathcal{F}'_{tf}) \geq \chi(\mathcal{O}_{Z_{1}}(a)) = a+1.
\end{equation*}
If $\mathcal{F}'_{tf}=0$ we are done.  Otherwise, we can repeat the argument on $\mathcal{F}'_{tf}$.

\item $s=0$, $t_{1} = 0$, and $0 \leq t_{2} \leq 2$. This is analogous to the previous case.

\item $s=0$, $t_{1} = t_{2} = 1$.  Let $\mathcal{O}_{Z}(a,b)$ be a summand of $\mathcal{E}$ such that the map to $\mathcal{F}$ does not vanish.  There are two cases to consider.  First suppose that $g|_{\mathcal{O}_{Z}(a,b)}$ is injective at the generic point of both components of $Z$.  Then $g|_{\mathcal{O}(a,b)}$ is injective and its cokernel is torsion.  Thus $\chi(\mathcal{F}) \geq \chi(\mathcal{O}_{Z}(a,b)) = a+b+1$.

Next, suppose that $g|_{\mathcal{O}_{Z}(a,b)}$ fails to be injective at the generic point of some component of $Z$.  Without loss of generality we may suppose that $g|_{\mathcal{O}_{Z}(a,b)}$ vanishes along $Z_{2}$, so that it must factor through $\mathcal{O}_{Z_{1}}(a)$.  Let $\mathcal{F}'$ denote the cokernel of $\mathcal{O}_{Z}(a,b) \to \mathcal{F}$.  Then
\begin{equation*}
\chi(\mathcal{F}) - \chi(\mathcal{F}'_{tf}) \geq \chi(\mathcal{O}_{Z_{1}}(a)) = a+1.
\end{equation*}
We then apply the argument above to $\mathcal{F}'_{tf}$.

\item $s=1$ and $t_{1}=t_{2}=0$. There are two options.  First suppose there is a summand $\mathcal{O}_{Z}(a,b)$ of $\mathcal{E}$ such that the restriction of $g$ to this summand does not vanish along either component of $Z$.  Then $g|_{\mathcal{O}_{Z}(a,b)}$ is injective, so that
\begin{equation*}
\chi(\mathcal{F}) \geq \chi(\mathcal{O}_{Z}(a,b)) = a+b+1 = a+b+\rk(\mathcal{G}).
\end{equation*}
Otherwise, there are different summands $\mathcal{O}_{Z}(a_{i},b_{i})$ for $i=1,2$ such that $g|_{\mathcal{O}_{Z}(a_{i},b_{i})}$ does not vanish along the generic point of $Z_{i}$.  Arguing as before, we see that the map $\mathcal{O}_{Z}(a_{1},b_{1}) \oplus \mathcal{O}_{Z}(a_{2},b_{2}) \to \mathcal{F}$ factors through a map
\begin{equation*}
\mathcal{O}_{Z_{1}}(a_{1}) \oplus \mathcal{O}_{Z_{2}}(b_{2}) \to \mathcal{F}.
\end{equation*}
Note that for each summand on the left the map to $\mathcal{F}$ must vanish at the node.  In particular, we see that this map has a non-trivial cokernel. Thus
\begin{align*}
\chi(\mathcal{F}) & \geq (a_{1} +1) + (b_{2} + 1) + 1 \\
& = (a_{1} +1) + (b_{2} + 1) + \rk(\mathcal{G}).
\end{align*}
\end{enumerate}

We next prove the induction step.  Choose any direct summand $\mathcal{O}_{Z}(a,b)$ of $\mathcal{E}$ such that the map $g: \mathcal{O}_{Z}(a,b) \to \mathcal{F}$ does not vanish at the node $p$.  Let $\mathcal{F}'$ denote the quotient of $\mathcal{F}$ by the $g$-image of $\mathcal{O}_{Z}(a,b)$.  We write
\begin{equation*}
\mathcal{F}'_{tf} = \mathcal{G}' \oplus \mathcal{H}'_{1} \oplus \mathcal{H}'_{2}
\end{equation*}
and $s',t'_{1},t'_{2}$ for the ranks.  Then one of the following must hold:

\begin{enumerate}
\item The map $g: \mathcal{O}_{Z}(a,b) \to \mathcal{F}$ is injective.  We deduce that
\begin{align*}
\chi(\mathcal{F}) - \chi(\mathcal{F}'_{tf}) & = (\chi(\mathcal{F}) - \chi(\mathcal{F}')) + (\chi(\mathcal{F}') - \chi(\mathcal{F}'_{tf})) \\
& \geq (a+b+1) + \chi(\mathcal{F}'_{tor}|_{p}).
\end{align*}
By Lemma \ref{lemm:gnonvanishing}, we see that
\begin{equation*}
\chi(\mathcal{F}) - \chi(\mathcal{F}'_{tf}) \geq (a+b) + (\rk(\mathcal{G}) - \rk(\mathcal{G}')).
\end{equation*}
Let $\mathcal{E}'$ denote the quotient of $\mathcal{E}$ by the chosen direct summand $\mathcal{O}_{Z}(a,b)$.  We have a surjection $\mathcal{E}' \to \mathcal{F}'_{tf}$.  By combining the equation above with the induction hypothesis for the surjection $\mathcal{E}' \to \mathcal{F}'_{tf}$ we deduce the desired inequality for $\chi(F)$.  By calculating ranks along the generic points of $Z_{1}$ and $Z_{2}$ we see that
\begin{equation*}
s + t_{1} = 1 + s' + t'_{1} \qquad \qquad s+t_{2} = 1 + s' + t'_{2}.
\end{equation*}
and the two rank equations can be deduced by the induction hypothesis.

\item The map $g: \mathcal{O}_{Z}(a,b) \to \mathcal{F}$ is not injective but the restriction to the generic point of $Z_{1}$ is injective.  Then $g$ must factor through $\mathcal{O}_{Z_{1}}(a)$.  Thus
\begin{align*}
\chi(\mathcal{F}) - \chi(\mathcal{F}'_{tf}) & = (\chi(\mathcal{F}) - \chi(\mathcal{F}')) + (\chi(\mathcal{F}') - \chi(\mathcal{F}'_{tf})) \\
& \geq (a+1) + \chi(\mathcal{F}'_{tor}|_{p}).
\end{align*}
By Lemma \ref{lemm:gnonvanishing}, we see that
\begin{equation*}
\chi(\mathcal{F}) - \chi(\mathcal{F}'_{tf}) \geq (a+1) + (\rk(\mathcal{G}) - \rk(\mathcal{G}')).
\end{equation*}
Let $\mathcal{E}'$ denote the quotient of $\mathcal{E}$ by the chosen direct summand $\mathcal{O}_{Z}(a,b)$.  We have a surjection $\mathcal{E}' \to \mathcal{F}'_{tf}$.  By combining the equation above with the induction hypothesis for the surjection $\mathcal{E}' \to \mathcal{F}'_{tf}$ we deduce the desired inequality for $\chi(F)$.  By calculating ranks along the generic points of $Z_{1}$ and $Z_{2}$ we see that
\begin{equation*}
s + t_{1} = 1 + s' + t'_{1} \qquad \qquad s+t_{2} = s' + t'_{2}.
\end{equation*}
and the two rank equations can be deduced by the induction hypothesis.

\item The map $g: \mathcal{O}_{Z}(a,b) \to \mathcal{F}$ is not injective but the restriction to the generic point of $Z_{2}$ is injective.  The argument is exactly analogous to the previous case.
\end{enumerate}
\end{proof}

\begin{example} \label{exam:pointblowup1}
Consider the vector bundle $\mathcal{E} = \mathcal{O}_{Z}(2,-1) \oplus \mathcal{O}_{Z}(-1,2)$ on $Z$.  Suppose that $\mathcal{F}$ is a torsion free sheaf with rank $1$ along each component of $Z$ such that there is a surjection $\mathcal{E} \to \mathcal{F}$.  There are two cases:
\begin{enumerate}
\item If $\mathcal{F}$ is locally free, then Lemma \ref{lemm:surjectiononnodalcurve} shows that $\chi(\mathcal{F}) \geq 2$.  This equality will be achieved only when $\mathcal{F}$ is one of the components of $\mathcal{E}$.
\item If $\mathcal{F}$ is not locally free, then Lemma \ref{lemm:surjectiononnodalcurve} shows that $\chi(\mathcal{F}) \geq 0$.  This equality will be achieved when $\mathcal{F} = \OO_{Z_{1}}(-1) \oplus \OO_{Z_{2}}(-1)$ and the surjection is constructed by taking the projection onto the $-1$-component on each factor.
\end{enumerate}
In both cases the bound given by Lemma \ref{lemm:surjectiononnodalcurve} is sharp.
\end{example}

In applications we are given a vector bundle $\mathcal{E}$ on a smoothing of $Z$ and would like to control the degrees along a general fiber of the sheaves $\mathcal{F}$ admitting a surjection $\mathcal{E} \to \mathcal{F}$.   We will show that the lowest possible degree of a torsion free sheaf $\mathcal{F}$ that has rank $m$ on a general fiber and admits a surjection $\mathcal{E} \to \mathcal{F}$ is bounded below by the following constant.

\begin{definition} \label{defi:degreebound}
Let $\mathcal{E}$ denote a vector bundle on the nodal curve $Z$ and let $\mathcal{E} = \oplus_{i=1}^{r} \mathcal{E}_{i}$ with $\mathcal{E}_{i} = \mathcal{O}(a_{i}, b_{i})$ be the decomposition into direct summands.  For any positive integer $m$, we define the degree bound
\begin{equation*}
\degbd(\mathcal{E},m) = \inf_{J,K_{1},K_{2}} \left\{  \sum_{i \in J}( a_{i} + b_{i} ) + \sum_{i \in K_{1}} (a_{i} + 1) + \sum_{i \in K_{2}} (b_{i} + 1) \right\}
\end{equation*}
where $J, K_{1}, K_{2}$ denote disjoint subsets of $\{ 1, 2, \ldots, r\}$ such that $|J| + |K_{1}| = m$ and $|J| + |K_{2}| = m$.
\end{definition}

The index set for the infimum represents the summands of $\mathcal{E}$ which could admit a generically surjective map to a sheaf $\mathcal{F}$ of rank $(m,m)$ as in Lemma \ref{lemm:surjectiononnodalcurve}.  For each summand the corresponding contribution to the infimum is similar to the bound given in Lemma \ref{lemm:surjectiononnodalcurve}, but Lemma \ref{lemm:restrictionoftfsheaf} will allow us to make a small improvement.

We will mostly be interested in the special case $m=1$, where we have
\begin{equation*}
\degbd(\mathcal{E},1) = \min\{ \inf_{i=1,\ldots,r}(a_{i}+b_{i}), \inf_{i=1,\ldots,r}(a_{i}) + \inf_{i=1,\ldots,r}(b_{j}) + 2 \}.
\end{equation*}

\begin{theorem} \label{theo:degreebound}
Let $\pi: \mathcal{C} \to B$ be a smoothing of $Z$ with smooth total family and let $\mathcal{E}$ be a vector bundle on $\mathcal{C}$.  Suppose that for a general curve $C$ in our family $\mathcal{E}|_{C}$ admits a surjection onto a vector bundle $\mathcal{Q}_{C}$ of rank $m$ and degree $d$.  Then $d \geq \degbd(\mathcal{E}|_{Z},m)$.
\end{theorem}

\begin{proof}
Let $\eta$ denote the generic point of $B$ and let $\mathcal{C}_{\eta}$ denote the generic fiber of $\pi$.  Since $\mathcal{C}_{\eta} \cong \mathbb{P}^{1}_{k(B)}$, the restriction of $\mathcal{E}$ to $\mathcal{C}_{\eta}$ decomposes into a direct sum of line bundles (in exactly the same way that $\mathcal{E}|_{C}$ decomposes).  Thus $\mathcal{E}|_{\mathcal{C}_{\eta}}$ admits a surjection onto a sheaf $\mathcal{Q}_{\eta}$ which has the same decomposition structure as $\mathcal{Q}_{C}$.  Spreading out, we find an open subset $U \subset B$ and a sheaf $\mathcal{Q}_{U}$ on $\pi^{-1}(U)$ such that there is a surjection $\mathcal{E}|_{\pi^{-1}(U)} \to \mathcal{Q}_{U}$.

Using properness of Quot schemes, we see that there is a torsion free sheaf $\mathcal{F}$ on $\mathcal{C}$ and a surjection $g: \mathcal{E} \to \mathcal{F}$ whose restriction to the general fiber $C$ coincides with the original surjection $\mathcal{E}|_{C} \to Q_{C}$.  Restricting to the central fiber $Z$, we still have a surjection
\begin{equation*}
\mathcal{E}|_{Z} \to \mathcal{F}|_{Z} \to (\mathcal{F}|_{Z})_{tf}.
\end{equation*}
Let $s$ denote the rank of the locally free part of $(\mathcal{F}|_{Z})_{tf}$ as in the statement of Lemma \ref{lemm:restrictionoftfsheaf}.  By Lemma \ref{lemm:restrictionoftfsheaf}, we see that $\chi((\mathcal{F}|_{Z})_{tf}) \leq d +s$.  By Lemma \ref{lemm:surjectiononnodalcurve} we find a direct summand $\oplus_{i=1}^{t} \mathcal{O}_{Z}(a_{i},b_{i})$ of $\mathcal{E}|_{Z}$ such that
\begin{equation*}
\sum_{i=1}^{j} (a_{i} + b_{i}) + \sum_{i=j+1}^{j+k_{1}} (a_{i}+1) + \sum_{i=j+k_{1}+1}^{j+k_{1}+k_{2}} (b_{i}+1) \leq \chi((\mathcal{F}|_{Z})_{tf}) - s.
\end{equation*}
Combining everything, we see that
\begin{align*}
\degbd(\mathcal{E}|_{Z},m) & \leq  \sum_{i=1}^{j} (a_{i} + b_{i}) + \sum_{i=j+1}^{j+k_{1}} (a_{i}+1) + \sum_{i=j+k_{1}+1}^{j+k_{1}+k_{2}} (b_{i}+1) \\
& \leq \chi((\mathcal{F}|_{Z})_{tf}) - s  \\
& \leq d+ s - s \\
& = d
\end{align*}
proving the desired statement.
\end{proof}

\begin{example}
Let $X$ be the blow-up of $\mathbb{P}^{1} \times \mathbb{P}^{1}$ at a point.  Consider the fibers of a projection map $f: X \to \mathbb{P}^{1}$; there is a single reducible fiber $Z$ which is the union of two rational curves at a single node.  In this case $T_{X}|_{Z} \cong \mathcal{O}_{Z}(2,-1) \oplus \mathcal{O}_{Z}(-1,2)$.   Let $C$ denote a general fiber of $f$.  Theorem \ref{theo:degreebound} predicts that every one-dimensional quotient of $T_{X}|_{C}$ will have non-negative degree.    This is of course easy to verify directly since $T_{X}|_{C} \cong \mathcal{O}_{C}(2) \oplus \mathcal{O}_{C}$.

It is interesting to observe explicitly how the argument of Theorem \ref{theo:degreebound} applies in this situation.  The infimum in Theorem \ref{theo:degreebound} is obtained by taking the $(-1)$-component in each summand of $T_{X}|_{Z}$.  Comparing against Example \ref{exam:pointblowup1}, we should expect to observe the quotient $\OO_{Z_1}(-1) \oplus \OO_{Z_2}(-1)$ as a limit of quotients on a general fiber.  Indeed, the map $T_{X}|_{C} \to f^{*}T_{\mathbb{P}^{1}}$ fails to be an isomorphism at the node in $Z$, yielding a surjection $T_{X}|_{C} \to \mathcal{I}_{p} \otimes f^{*}T_{\mathbb{P}^{1}}$.  The restriction of this map to a general fiber is $T_{X}|_{C} \to \mathcal{O}_{C}$; the restriction to the special fiber is a map $T_{X}|_{Z} \to \mathcal{F}$ where $\mathcal{F}$ has torsion part of length $1$ and has torsion-free part $\OO_{Z_1}(-1) \oplus \OO_{Z_2}(-1)$.
\end{example}

We next show that Theorem \ref{theo:degreebound} is in some sense the optimal result possible.  The following proposition identifies the key construction.

\begin{proposition} \label{clai:rank2construction}
Let $S$ denote the blow-up of $\mathbb{P}^{1} \times \mathbb{P}^{1}$ at a point equipped with one of the projection maps $\pi: S \to \mathbb{P}^{1}$.  Let $Z$ denote the nodal fiber with components $Z_{1},Z_{2}$ and let $p$ denote the nodal point. Let $T_{1}$ be a section of $\pi$ that meets $Z_{1}$ and let $T_{2}$ be a section of $\pi$ that meets $Z_{2}$.
Choose integers $a,b,c,d$ satisfying $c \geq a+2$ and $b \geq d+2$.  There is a rank $2$ vector bundle $\mathcal{E}_{a,b,c,d}$ on $S$ that fits into an exact sequence
\begin{equation*}
0 \to \mathcal{O}_{S}((c-1)T_{1} + (b-1)T_{2}) \to \mathcal{E}_{a,b,c,d} \xrightarrow{\psi} \mathcal{I}_{p}((a+1)T_{1} + (d+1)T_{2}) \to 0
\end{equation*}
so that
\begin{enumerate}
\item The restriction  of the map $\psi$ to a general fiber $C$ of $\pi$ yields the surjection $\mathcal{O}_{C}(b+c-2) \oplus \mathcal{O}_{C}(a+d+2) \to \mathcal{O}_{C}(a+d+2)$.
\item The restriction of the map $\psi$ to $Z$ yields the surjection $\mathcal{O}_{Z}(a,b) \oplus \mathcal{O}_{Z}(c,d) \to \mathcal{O}_{Z_{1}}(a) \oplus \mathcal{O}_{Z_{2}}(d)$.
\end{enumerate}
\end{proposition}

\begin{proof}
We first recall the Serre construction. If $\mathcal{L}$ is a line bundle on $S$ satisfying $H^{2}(S,\mathcal{L}) = 0$, then Serre's construction (see for example the summary at the end of \cite{Schnell}, with $Z = p$, $X = S$) yields a locally free sheaf $\mathcal{E}$ obtained via an extension
\begin{equation*}
0 \to \mathcal{L} \to \mathcal{E} \to \mathcal{I}_{p} \to 0.
\end{equation*}

If $D$ is an effective divisor, then $h^{2}(S,\mathcal{O}_{S}(D-Z)) \leq h^{2}(S,\mathcal{O}_{S}(-Z)) = 0$. Thus, we may perform the Serre construction above with $\mathcal{L} = \OO_S((c-a-2)T_1+(b-d-2)T_2-Z)$, and twist by $\OO_{S}((a+1)T_1+(d+1)T_2)$ to obtain the sequence
\begin{equation*}
0 \to \mathcal{O}_{S}((c-1)T_{1} + (b-1)T_{2}) \to \mathcal{E}_{a,b,c,d} \to \mathcal{I}_{p}((a+1)T_{1} + (d+1)T_{2}) \to 0.
\end{equation*}

The restriction of the above sequence to $C$ is
\[ 0 \to \OO_C(b+c-2) \to \mathcal{E}_{a,b,c,d} \to \OO_C(a+d+2) \to 0 \]
so it only remains to calculate the restriction to $Z$.  Restricting to $Z_{1}$ yields
\begin{equation*}
0 \to \mathcal{O}_{Z_{1}}(c-1) \to \mathcal{E}_{a,b,c,d}|_{Z_{1}} \to \mathcal{O}_{Z_{1}}(a) \oplus k(p) \to 0
\end{equation*}
and since by assumption $c \geq a+2$ we see that $\mathcal{E}|_{Z_{1}} \cong \mathcal{O}_{Z_{1}}(c) \oplus \mathcal{O}_{Z_{1}}(a)$.  Similarly $\mathcal{E}_{a,b,c,d}|_{Z_{2}} = \OO_{Z_2}(b) \oplus \OO_{Z_2}(d)$. This means that the sequence restricted to $Z$ is
\begin{equation}  \label{eq:serreconst}
0 \to \OO_{Z}(c-1,b-1) \to \mathcal{E}_{a,b,c,d}|_Z \to \OO_{Z_1}(a) \oplus \OO_{Z_2}(d) \to 0.
\end{equation}
Note that the $\OO_{Z_{1}}(c)$ summand of $\mathcal{E}_{a,b,c,d}|_{Z_{1}}$ and the $\OO_{Z_{2}}(b)$ summand of $\mathcal{E}_{a,b,c,d}|_{Z_{2}}$ cannot line up since if they did the rightmost term in Equation \eqref{eq:serreconst} would need to have a torsion subsheaf supported in dimension $0$.
Altogether we see that
\begin{equation*}
\mathcal{E}|_{Z} \cong \mathcal{O}_{Z}(a,b) \oplus \mathcal{O}_{Z}(c,d).
\end{equation*}
\end{proof}

The next theorem shows that the degree bound given in Theorem \ref{theo:degreebound} is sharp.

\begin{theorem}
Let $\pi: S \to \mathbb{P}^{1}$ be the composition of the blow-up of a point $\phi: \mathcal{C} \to \mathbb{P}^{1} \times \mathbb{P}^{1}$ with a projection map.  Let $Z$ be the central fiber with components $Z_{1},Z_{2}$ and node $p$ and fix a vector bundle $\mathcal{E}_{Z}$ on $Z$.  For any $m$ satisfying $0 \leq m \leq \rk(\mathcal{E}_{Z})$ there is
\begin{itemize}
\item a vector bundle $\mathcal{E}$ on $S$ of rank $\rk(\mathcal{E}|_{Z})$ such that $\mathcal{E}|_{Z} = \mathcal{E}_{Z}$, and
\item a reflexive sheaf $\mathcal{F}$ on $S$ whose restriction to a general fiber $C$ has rank $m$ and degree $\degbd(\mathcal{E}|_{Z},m)$
\end{itemize}
such that $\mathcal{E}$ admits a surjection onto $\mathcal{F}$.
\end{theorem}

\begin{proof}
We let $E$ denote the exceptional divisor of $\phi$, let $H_{1}$ denote a divisor representing $\pi^{*}\mathcal{O}_{\mathbb{P}^{1}}(1)$, and let $H_{2}$ denote a general fiber of the other projection map to $\mathbb{P}^{1}$. We let $T_{1}$ be a section of $\pi$ that meets $Z_{1}$ and let $T_{2}$ be a section of $\pi$ that meets $Z_{2}$.

Let $\mathcal{E}'_{Z} \subset \mathcal{E}_{Z}$ be a direct summand which realizes the infimum in $\degbd(\mathcal{E}_{Z},m)$.  It suffices to prove the statement for $\mathcal{E}'_{Z}$.  We write
\begin{equation*}
\mathcal{E}'_{Z} = \bigoplus_{i=1}^{j} \mathcal{O}_{Z}(a_{i},b_{i}) \oplus \bigoplus_{i = j+1}^{j+k} \mathcal{O}_{Z}(a_{i},b_{i}) \oplus \bigoplus_{i = j+k+1}^{j+2k} \mathcal{O}_{Z}(a_{i},b_{i})
\end{equation*}
where $j+k = m$ and the decomposition coheres to the notation used in Definition \ref{defi:degreebound}.  Since we are assuming that $\mathcal{E}'$ computes $\degbd(\mathcal{E}_{Z},m)$ we must have that $b_{i} \geq b_{i+k}+2$ and $a_{i}+2 \leq a_{i+k}$ whenever $j+1 \leq i \leq j+k$.

We construct the desired global sheaf $\mathcal{E}'$ and the surjection $\mathcal{E}' \to \mathcal{F}$ using direct sums.
\begin{itemize}
\item For $1 \leq i \leq j$, let $\mathcal{E}_{i}$ be the line bundle $\mathcal{O}_{S}((a_{i}+b_{i})H_{2} - a_{i}E)$.  Set $\mathcal{F}_{i} = \mathcal{E}_{i}$.  Let $\psi_{i}: \mathcal{E}_{i} \to \mathcal{F}_{i}$ be the identity map.
\item For $j+1 \leq i \leq j+ k$, let $\mathcal{E}_{i}$ be the rank $2$ bundle $\mathcal{E}_{a_{i},b_{i},a_{i+k},b_{i+k}}$ constructed in Claim \ref{clai:rank2construction}.  (Note that we have already verified the required inequalities on the indices.)  Setting $\mathcal{F}_{i} = \mathcal{I}_{p}((a_{i}+1)T_{1} + (b_{i+k}+1)T_{2})$, the claim shows that we have a surjection $\mathcal{E}_{i} \to \mathcal{F}_{i}$.
\end{itemize}
Set $\mathcal{F} = \oplus_{i=1}^{j+k} \mathcal{F}_{i}$.  It only remains to show that the degree of $\mathcal{F}|_{C}$ is $\degbd(\mathcal{E}|_{Z},m)$.  For $1 \leq i \leq j$ we have that $\mathcal{F}_{i}$ contributes $a_{i} + b_{i}$ to the degree, and for $j+1 \leq i \leq j+k$ we have that $\mathcal{F}_{i}$ contributes $a_{i} + b_{i+k} + 2$ to the degree, finishing the proof.
\end{proof}

\subsection{Gluing results}

We end this section by applying these results to control the behavior of a locally free sheaf for a smoothing of a union of two free curves.
As we discussed before, the key observation is that we can ensure that the total space of the deformation is smooth.

\begin{theorem} \label{theo:surjectionsandsmoothing}
Let $X$ be a smooth projective variety and let $\mathcal{E}$ denote a torsion-free sheaf on $X$.  Suppose $f_{1}: Z_{1} \to X$ and $f_{2}: Z_{2} \to X$ are two free rational curves in $X$ meeting at a point $p$ whose images are contained in the locally free locus of $\mathcal{E}$.  Let $f_{Z}: Z \to X$ denote the map obtained by gluing $Z_{1}$ and $Z_{2}$ at $p$ and let $f_{C}: C \to X$ be a general smoothing of $Z$.  For any integer $m$ satisfying $1 \leq m \leq \dim X$ the vector bundle $f_{C}^{*}\mathcal{E}$ does not have any summands with rank $m$ and degree less than $\degbd(f_{Z}^{*}\mathcal{E},m)$.
\end{theorem}

\begin{proof}
By Lemma \ref{lemm:smoothingforfreecurves} there is a smooth surface $S$ equipped with a flat projective morphism to a curve $\pi: S \to T$, a point $0 \in T$, and a morphism $g: S \to X$ such that the restriction of $g$ to the central fiber is $f_{Z}$.  We apply Theorem \ref{theo:degreebound} to $g^{*}\mathcal{E}$.
\end{proof}

\section{Stability for low rank vector bundles}
In this section we explore the relationship between the stability of vector bundles and the stability of their restrictions to free rational curves.  Recall that we can discuss stability with respect to any nef curve class $\alpha$ by using the slope $\mu_{\alpha}(\mathcal{E}) = \frac{c_1(\mathcal{E}) \cdot \alpha}{\rk(\mathcal{E})}$.  We start with a familiar and elementary observation.

\begin{proposition} \label{prop:restrictssthenbundless}
Let $X$ be a smooth projective variety and let $\mathcal{E}$ be a vector bundle on $X$.   Suppose that a variety $M$ parametrizes a dominant family of irreducible curves on $X$ such that for any codimension $2$ closed subset $Z \subset X$ there is a curve parametrized by $M$ that is disjoint from $Z$.  Let $C$ be a general curve parametrized by $M$ and let $\alpha$ denote its numerical class.  If $\mathcal{E}|_C$ is (semi-)stable, then $\mathcal{E}$ is (semi-)stable on $X$ with respect to $\alpha$. 
\end{proposition}

\begin{proof}
We prove the contrapositive statement.  Suppose $\mathcal{E} \to \mathcal{Q}$ is a destabilizing quotient for $\mathcal{E}$ with respect to $\alpha$. Then restricting $\mathcal{E}$ and $\mathcal{Q}$ to $C$ gives a surjection $\mathcal{E}|_C \to \mathcal{Q}|_C$. Computing slopes we get
\[ \mu(\mathcal{E}|_C) = \frac{c_1(\mathcal{E}|_C)}{\rk \mathcal{E}|_C} = \frac{c_1(\mathcal{E}) \cdot \alpha}{\rk \mathcal{E}} = \mu_{\alpha}(\mathcal{E}). \]
By generality $C$ is contained in the locally free locus of $\mathcal{Q}$, and in particular $\mathcal{Q}|_{C}$ is torsion free.  Thus a similar argument shows that $\mu(\mathcal{Q}|_{C}) =  \mu_{\alpha}(\mathcal{Q})$.  Thus $\mathcal{E}|_C \to \mathcal{Q}|_C$ is also destabilizing.  The argument for semi-stability is the same.
\end{proof}

\begin{corollary} \label{coro:0balancedisss}
If $X$ has a $0$-balanced free rational curve $C$, then $T_X$ must be semistable with respect to the numerical class of $C$.
\end{corollary}

\begin{proof}
Since the $0$-balanced property is preserved under generalization, this follows from Proposition \ref{prop:restrictssthenbundless}.
\end{proof}

We will be interested in identifying situations where a converse statement holds.  This question was also analyzed by \cite{Ran01} when the ambient variety is $\mathbb{P}^{n}$.  We will need the following (well-known) result concerning $k$-planes in $\PP^{n}$.

\begin{lemma} \label{lemm:kplanelemma}
Consider an irreducible family of $k$-planes in $\PP^{n}$.  Suppose that the intersection of any two planes in our family has dimension $(k-1)$.  Then either all the $k$-planes in our family are contained in a fixed $(k+1)$-dimensional subspace of $\PP^{n}$ or there is a fixed $(k-1)$-dimensional subspace of $\PP^{n}$ which is contained in every $k$-plane.
\end{lemma}

\begin{proof}
Suppose that there is no fixed $(k-1)$ plane that is contained in every $k$-plane in our family.  Fix $2$ general planes $L_{1},L_{2}$ in our family.  Since a general $k$-plane $L$ in our family does not contain $L_{1} \cap L_{2}$, it must be contained in the span of $L_{1}$ and $L_{2}$ which is a $(k+1)$-plane.
\end{proof}

The following result is the main theorem in this section.  As discussed in the introduction, it can be seen as a ``Mehta-Ramanathan type statement'' for families of free rational curves and torsion free sheaves of low rank.

\begin{theorem} \label{theo:balanceduptorank6}
Let $X$ be a smooth projective variety and let $\mathcal{E}$ be a torsion free sheaf on $X$ of rank $r \leq 5$.  Suppose that $\alpha$ is the class of a free rational curve $C_{0}$ that is general in an irreducible component of $\Rat_{conn}(X)$.  There is some constant $\Gamma = \Gamma(r)$ such that $\mathcal{E}$ is semistable with respect to $\alpha$ if and only if there is a curve $\widetilde{C}$ obtained by gluing and smoothing $\Gamma(r)$ copies of $C_{0}$ with $\mathcal{E}|_{\widetilde{C}}$ is semistable.
\end{theorem}

\begin{remark}
Loosely speaking, our strategy in the proof of Theorem \ref{theo:balanceduptorank6} is to think about whether the positive directions in $\mathcal{E}|_C$ vary a lot as we deform $C$. If they vary a lot, we can glue two curves of class $C$ together and deform to a curve $C'$ with more balanced $\mathcal{E}|_{C'}$. If they don't vary, then we argue that they correspond to a positive subsheaf of $\mathcal{E}$.

However, one cannot expect a precise relationship along these lines even when $\mathcal{E}$ is the tangent bundle of $X$. For instance, let $X$ be a homogeneous contact manifold coming from an orthogonal or exceptional simple Lie algebra as described in Section 1.4 of \cite{HwangMok}. Then $X$ has dimension $2m+1$ and has a single minimal family of rational curves, which are lines under the natural projective embedding.  Let $C$ be a general such line. Then $T_X|_C = \OO(2) \oplus \OO(1)^m \oplus \OO^m$, with the $\OO(2)$ factor corresponding to the image of $T_C$ in $T_X|_C$. At a point $x$ of $X$, the tangent spaces to the lines span a codimension one subspace $D_x \subset T_{X,x}$, giving rise to a distribution $D \subset T_X$. Even though $D$ is the span of the ``positive'' directions of these rational curves, we have that $D|_C = \OO(2) \oplus \OO(1)^{m-1} \oplus \OO^{m-1} \oplus \OO(-1)$, which is not destabilizing for $T_X$. Thus, even the distribution spanned by the positive directions of a family of curves can be less ``positive'' than the tangent bundle.
\end{remark}

Before starting the proof let us give a couple reminders of facts established earlier in the paper.  First, the restriction of a torsion free sheaf to a general member of a family of free curves will be locally free, and we will use this fact without further mention.  Second, when we speak of ``gluing and smoothing'' free rational curves $C_{1},C_{2}$ which meet at a point we allow the case $C_{1} = C_{2}$ (in which case the resulting stable map could be a double cover of the underlying curve). 

\begin{proof}
Since the reverse implication is established by Corollary \ref{coro:0balancedisss}, we only need to prove the forward implication.  The proof is by induction on $r$.  The base case when $r = 1$ is trivial.

When $r > 1$, we begin by making some reductions.  First, observe that by Lemma \ref{lem-glueNonACovers} we can replace $C_{0}$ by a general smoothing of a gluing of $r$ curves which are deformations of $C_{0}$ (and absorb the factor of $r$ into $\Gamma(r)$).  This ensures that the slope $\mu(\mathcal{E}|_{C_0})$ is an integer which we denote by $a$.  Furthermore, by replacing $C_{0}$ by gluing and smoothing $\mathrm{lcm}(\Gamma(2),\ldots,\Gamma(r-1))$ copies of $C_{0}$ (and absorbing the factor into $\Gamma(r)$) we may assume that any semistable sheaf $\FF$ of rank $< r$ has the property that $\FF|_{C_{0}}$ is semistable.

We will always assume that $C_{0}$ is general in its family, so that the direct summand decomposition of $\mathcal{E}|_{C_{0}}$ will be the same as the direct summand decomposition of the restriction of $\mathcal{E}$ to a general deformation of $C_{0}$. The version of the Grauert-Mulich Theorem proved in \cite[Theorem 3.1]{PatelRiedlTseng} shows that $\mathcal{E}|_{C_{0}}$ will be sequential.

\bigskip

$\mathbf{r=2}$: By Lemma \ref{lemm:citedGrauert-Mulich} the only possibility for $\mathcal{E}|_{C_0}$ is $\OO(a) \oplus \OO(a)$ which is semistable.

\bigskip

$\mathbf{r=3}$: By Lemma \ref{lemm:citedGrauert-Mulich} $\mathcal{E}|_{C_0}$ is isomorphic to either $\OO(a)^3$ or $\OO(a+1) \oplus \OO(a) \oplus \OO(a-1)$, and we only need to consider the latter case.  Fix a general point $p$ of $X$ and consider the sublocus $S_{p} \subset \mathcal{E}|_{p}$ swept out by the $\OO(a+1)$ summands as we vary $C$ over all deformations of $C_{0}$ through $p$ such that $\mathcal{E}|_{C} \cong \OO(a+1) \oplus \OO(a) \oplus \OO(a-1)$.  We claim that $S_{p}$ spans all of $\mathcal{E}|_p$.  Suppose for a contradiction that this is not the case.  Let $\FF$ be the saturated subsheaf of $\mathcal{E}$ constructed in Lemma \ref{lemm:positiveSubsheaf} whose restriction to a general point $p$ is the span of $S_{p}$.  Letting $C$ denote a general deformation of $C_{0}$, Lemma \ref{lemm:positiveSubsheaf} shows that $\FF|_{C}$ will contain the $\OO(a+1)$ factor of $\mathcal{E}|_{C}$.  Our induction assumption implies that the Harder-Narasimhan filtration of $\FF$ restricts to give the Harder-Narasimhan filtration of $\FF|_{C}$.  Since $\FF|_{C}$ has at least one summand of degree $\geq a+1$, we see that the maximal destabilizing subsheaf of $\FF$ will also destabilize $\mathcal{E}$.  This verifies the claim.

Suppose that $C_{1}$ and $C_{2}$ are two general curves in the deformation class of $C_0$ that pass through the same general point $p$.  The argument above shows that the $\OO(a+1)$ summand in $C_{1}$ will line up with the $\OO(a-1)$ summand in $C_{2}$ (and vice versa).  By Theorem \ref{theo:surjectionsandsmoothing} we deduce that a general smoothing $\widetilde{C}$ of the gluing of $C_{1}$ and $C_{2}$ will be $0$-balanced.

\bigskip

$\mathbf{r=4}$:  By Lemma \ref{lemm:citedGrauert-Mulich} $\mathcal{E}|_{C_0}$ is isomorphic to either $\OO(a)^4$ or $\OO(a+1) \oplus \OO(a)^2 \oplus \OO(a-1)$, and we need only consider the latter case.  Let $C_{1}$ and $C_{2}$ be two general deformations of $C_{0}$ through a fixed point $p$.  As in the $r=3$ case, the $\OO(a+1)$ factor on $C_{1}$ and the $\OO(a-1)$ factor on $C_{2}$ will line up unless there is a saturated subsheaf $\FF \subsetneq \mathcal{E}$ such that $\FF|_p$ is the span of the $\OO(a+1)$ directions of deformations of $C_{0}$ through a general point $p$.

Suppose for a contradiction that there is such a subsheaf $\FF$.  Note that for a general deformation $C$ of $C_{0}$ the restriction $\FF|_{C}$ contains the $\OO(a+1)$ factor in $\mathcal{E}|_{C}$.  By the induction assumption the Harder-Narasimhan filtration of $\FF$ will restrict to give the Harder-Narasimhan filtration of $\FF|_{C}$ and every term of the Harder-Narasimhan filtration of $\FF|_{C}$ will be $0$-balanced.  Since $\FF|_{C}$ has at least one summand of degree $\geq a+1$, we see that the maximal destabilizing subsheaf of $\FF$ will also destabilize $\mathcal{E}$.

\bigskip

$\mathbf{r=5}$: Aside from the semistable case, Lemma \ref{lemm:citedGrauert-Mulich} gives three possibilities for $\mathcal{E}|_C$.

\textbf{Case 1:} $\mathcal{E}|_C = \OO(a+1) \oplus \OO(a)^3 \oplus \OO(a-1)$.  Suppose that as we vary $C$ over general deformations of $C_{0}$ through a general point $p$ the $\OO(a+1)$ factors at $p$ span all of $\mathcal{E}|_{p}$.  Then we can argue as before.

Otherwise, let $\FF \subsetneq \mathcal{E}$ be the saturated subsheaf such that for a general point $p$ the restriction $\FF|_{p}$ is the span of the $\OO(a+1)$ summands for deformations of $C_{0}$.  Lemma \ref{lemm:positiveSubsheaf} guarantees that there is an inclusion $\mathcal{O}(a+1) \to \FF|_{C}$.  By induction we know that the restriction of the maximal destabilizing subsheaf of $\FF$ to $C$ is balanced.  Thus it will destabilize $\mathcal{E}$.

\textbf{Case 2:} $\mathcal{E}|_C = \OO(a+1)^2 \oplus \OO(a) \oplus \OO(a-1)^2$.  Let $\FF \subset \mathcal{E}$ denote the saturated subsheaf such that for a general point $p$ the restriction $\FF|_{p}$ is the span of the $\OO(a+1)^{2}$ subspaces at $p$ for general deformations $C$ of $C_{0}$.  By Lemma \ref{lemm:positiveSubsheaf} there is an inclusion $\OO(a+1)^{2} \to \FF|_{C}$ for a general deformation $C$ of $C_{0}$.  If the rank of $\FF$ is smaller than $5$, then by the induction assumption the Harder-Narasimhan filtration of $\FF$ will restrict to give the Harder-Narasimhan filtration of $\FF|_{C}$ and every term of the Harder-Narasimhan filtration of $\FF|_{C}$ will be $0$-balanced.  In this case the maximal destabilizing subsheaf of $\FF$ would also destabilize $\mathcal{E}$, a contradiction.  We conclude that the $\OO(a+1)^{2}$ subspaces sweep out all of $\mathcal{E}|_{p}$ for a general point $p$.

We also claim that for two general deformations $C_{1},C_{2}$ of $C$ through a general point $p$, the $\OO(a+1)^{2}$ subspaces for $C_{1}$ and $C_{2}$ at $p$ only intersect at $0$.  If not, then they must generically intersect in a $1$-dimensional subspace (since by the argument above as we vary $C$ their span is all of $\mathcal{E}|_{p}$).  By Lemma \ref{lemm:kplanelemma} we conclude that either there is a fixed three-dimensional subspace of $\mathcal{E}|_{p}$ which contains every $\OO(a+1)^{2}$ subspace or there is a fixed one-dimensional subspace of $\mathcal{E}|_{p}$ which is contained in every $\OO(a+1)^{2}$ summand.  Furthermore, the previous paragraph shows that the first case cannot happen.  Thus there is a distinguished one-dimensional subspace at every general point $p$.  As we vary $p$ these $1$-dimensional spaces yield a rank $1$ saturated subsheaf $\mathcal{G} \subset \mathcal{E}$.

Note that if we restrict $\mathcal{G} \to \mathcal{E}$ to a general deformation $C$ of $C_{0}$, then $\mathcal{G}|_{C}$ must map to the $\OO(a+1)^{2}$ summand in $\mathcal{E}|_{C}$.  Consider the exact sequence
\begin{equation*}
0 \to \mathcal{G} \to \mathcal{E} \to \mathcal{Q} \to 0.
\end{equation*}
Then one of the following two situations must occur:
\begin{enumerate}
\item $\mathcal{G}|_{C}$ has slope $a+1$.  Then it is a destabilizing subsheaf, a contradiction.
\item $\mathcal{G}|_{C}$ has slope $<a+1$.  Then $\mathcal{Q}|_{C} \cong \mathcal{O}(b) \oplus \OO(a) \oplus \OO(a-1)^{2}$ for some $b \geq a+2$.  By Lemma \ref{lemm:citedGrauert-Mulich} $\mathcal{Q}$ has a quotient sheaf $\mathcal{Q}'$ such that $\mathcal{Q}'|_{C} \cong \OO(a) \oplus \OO(a-1)^{2}$.  But then the quotient $\mathcal{E} \to \mathcal{Q}'$ is destabilizing, a contradiction.
\end{enumerate}
This completes the argument showing that the $\OO(a+1)^{2}$ subspaces only intersect at $0$.

Finally, we claim that the $\OO(a+1)^{2}$ subspace for $C_{1}$ at $p$ meets the $\OO(a+1)^{2} \oplus \OO(a)$ subspace for $C_{2}$ at $p$ transversally.  If not, then since the two $\OO(a+1)^{2}$ subspaces intersect only at $0$ the two $\OO(a+1)^{2} \oplus \OO(a)$ subspaces will intersect in dimension exactly $2$.  By Lemma \ref{lemm:kplanelemma} the $\OO(a+1)^{2} \oplus \OO(a)$ factors will either always be contained in a fixed four-dimensional subspace or will all intersect along a fixed two-dimensional subspace.  The first case cannot happen since the $\OO(a+1)^{2}$ subspaces span all of $\mathcal{E}|_{p}$.  If the second happens, then arguing as before one can deduce that $\mathcal{E}$ is not stable, a contradiction.

Altogether, suppose that $C_{1},C_{2}$ are two general deformations of $C$ through a fixed general point $p$.  Let $Z$ denote their union.  Then
\begin{equation*}
\mathcal{E}|_{Z} \cong \mathcal{O}_{Z}(a+1,a-1)^{2} \oplus \mathcal{O}_{Z}(a,a) \oplus \mathcal{O}_{Z}(a-1,a+1)^{2}.
\end{equation*}
By Theorem \ref{theo:surjectionsandsmoothing} we deduce that for a general smoothing $\widetilde{C}$ of $Z$ the restriction $\mathcal{E}|_{\widetilde{C}}$ will not admit a surjection to $\mathcal{O}(a-1)$.  We conclude that $\mathcal{E}|_{\widetilde{C}}$ is $0$-balanced. 

\textbf{Case 3:} $\mathcal{E}|_C = \OO(a+2) \oplus \OO(a+1) \oplus \OO(a) \oplus \OO(a-1) \oplus \OO(a-2)$.  By arguing as above, we see that as we vary $C$ over all general deformations of $C_{0}$ through a general point $p$ the $\OO(a+2)$ factors must sweep out all of $\mathcal{E}|_{p}$.  

We next claim that for two general deformations $C_{1},C_{2}$ of $C_{0}$ through a general point $p$, the $\OO(a+2) \oplus \OO(a+1)$ subspaces at $p$ will meet transversally.  If not, then by applying Lemma \ref{lemm:kplanelemma} as before we would obtain a rank $1$ saturated subsheaf $\mathcal{G}$ such that at a general point $p$ the restriction $\mathcal{G}|_{p}$ is contained in every $\OO(a+2) \oplus \OO(a+1)$ subspace.  But this would yield a destabilizing subsheaf or quotient of $\mathcal{E}$.

We next claim that for two general deformations $C_{1},C_{2}$ of $C_{0}$ through a general point $p$, the $\OO(a+2) \oplus \OO(a+1)$ subspace on $C_{1}$ at $p$ will meet the $\OO(a+2) \oplus \OO(a+1) \oplus \OO(a)$ subspace transversally.  If not, then the two $\OO(a+2) \oplus \OO(a+1) \oplus \OO(a)$ subspaces for the two curves would intersect along a two-dimensional subspace.  By combining the previous paragraph with Lemma \ref{lemm:kplanelemma} we see there is a fixed $2$-dimensional subspace of $\mathcal{E}|_{p}$ that is contained in every $\OO(a+2) \oplus \OO(a+1) \oplus \OO(a)$ subspace.  Again this would yield a destabilizing subsheaf or quotient.

Altogether, suppose that $C_{1},C_{2}$ are two general deformations of $C$ through a fixed general point $p$.  Let $Z$ denote their union.  Then
\begin{equation*}
\mathcal{E}|_{Z} \cong \mathcal{O}_{Z}(a+2,a-2) \oplus \mathcal{O}_{Z}(a+1,a-1) \oplus \mathcal{O}_{Z}(a,a) \oplus \mathcal{O}_{Z}(a-1,a+1) \oplus \mathcal{O}_{Z}(a-2,a+2).
\end{equation*}
By Theorem \ref{theo:surjectionsandsmoothing} we deduce that for a general smoothing $\widetilde{C}$ of $Z$ the restriction $\mathcal{E}|_{\widetilde{C}}$ will not admit a rank $2$ quotient of degree $\leq -3$.  By Lemma \ref{lemm:citedGrauert-Mulich}, if $\mathcal{E}|_{\widetilde{C}}$ admitted a surjection to $\mathcal{O}(a-2)$ then it would also need to admit a surjection to a rank $2$ bundle of degree $\leq -3$, which we just showed is not possible.  Thus every summand of $\mathcal{E}|_{\widetilde{C}}$ has degree at least $a-1$.  This means that the new curve $\widetilde{C}$ will either be in Case 1 or in Case 2 above, and we have reduced the statement to a known case.
\end{proof}

\begin{remark}
Our approach for Theorem \ref{theo:balanceduptorank6} cannot be extended to the case when $\mathcal{E}$ has rank $7$.  The issue occurs when we have a curve $C_{0}$ such that
\begin{equation*}
\mathcal{E}|_{C_{0}} \cong \mathcal{O}(a+3) \oplus \mathcal{O}(a+2) \oplus \mathcal{O}(a+1) \oplus \mathcal{O}(a) \oplus \mathcal{O}(a-1) \oplus \mathcal{O}(a-2) \oplus \mathcal{O}(a-3).
\end{equation*}
In the best case scenario, for two general deformations $C_{1},C_{2}$ through a point $p$ the $\mathcal{O}(a+b)$ summand in $C_{1}$ will line up with the $\mathcal{O}(a-b)$ summand in $C_{2}$.   Suppose we glue the two curves and smooth to get $\widetilde{C}$.  Even in this optimal situation, Theorem \ref{theo:surjectionsandsmoothing} does not prove that $\mathcal{E}|_{\widetilde{C}}$ is more balanced than $\mathcal{E}|_{C_{0}}$.  It only shows that $\mathcal{E}|_{\widetilde{C}}$ is not less balanced than $\mathcal{E}|_{C_{0}}$.
\end{remark}

\section{Factoring covers}

Our main tools for understanding direct summand decompositions of restricted bundles -- \cite[Proposition 3.1]{PatelRiedlTseng} and Theorem \ref{theo:balanceduptorank6} -- include the assumption that the evaluation map for a family of free rational curves has connected fibers.  In this section we analyze families of free rational curves for which the evaluation map fails to have connected fibers.  The key definition is the following.

\begin{definition} \label{defi:factoringmap}
Let $X$ be a smooth projective variety.  Let $M$ be an irreducible component of $\Rat(X)$.  Suppose that $\psi: Y \to X$ is a generically finite map from a projective variety $Y$ and that $N$ is an irreducible component of $\Rat(Y)$ such that composition with $\psi$ defines a dominant rational map $\psi_{*}: N \dashrightarrow M$.  Then we say that the component $N$ and the map $\psi: Y \to X$ form a \emph{factoring cover} for the component $M$.  (Often we will drop $N$ from the notation.)
\end{definition}

\begin{remark}
In the setting of Definition \ref{defi:factoringmap} it is possible that there are several different irreducible components of $\overline{M}_{0,0}(Y)$ which all map dominantly to the same irreducible component of $\overline{M}_{0,0}(X)$.   Thus when discussing factoring covers it is important to specify the component $N$ of $\Rat(Y)$.

For example, let $f: Y \to \Hilb^{2}(\mathbb{P}^{2})$ be the double cover obtained by blowing up the diagonal of $\mathbb{P}^{2} \times \mathbb{P}^{2}$.  For $i=1,2$ let $N_{i}$ denote the strict transform on $Y$ of the family of lines contained in the fibers of the $i$th projection map of $\mathbb{P}^{2} \times \mathbb{P}^{2}$.  Then $N_{1}$ and $N_{2}$ map to the same component $M$ of $\Rat(X)$ and both components give $Y$ the structure of a factoring cover.
\end{remark}

Our main example of factoring covers comes from irreducible components $M$ of $\Rat(X)$ such that the evaluation map $\mathcal{U} \to X$ for the one-pointed family $\mathcal{U}$ over $M$ fails to have connected fibers.  In this case the finite part of the Stein factorization of the evaluation map will be a factoring cover in a natural way.  However, there are many other examples of factoring covers.

\begin{example} \label{exam:birationalfactoringcover}
Let $X$ be a smooth projective variety and $M$ an irreducible component of $\Rat(X)$.  Suppose that $\psi: Y \to X$ is a factoring cover for $M$.  If we precompose $\psi$ by any birational map $\phi: Y' \to Y$ then the composition $\psi \circ \phi$ is still a factoring cover.  Indeed, we can simply take the strict transform on $Y'$ of the family of rational curves on $Y$ parametrized by $N$.
\end{example}

\begin{lemma} \label{lemm:tangentbundlefactoringcover}
Let $X$ be a smooth projective variety.  Suppose that $M$ is an irreducible component of $\Rat(X)$ and that $\psi: Y \to X$ is a factoring cover for $M$ equipped with an irreducible component $N$ of $\Rat(Y)$.  For a general map $g$ parametrized by $N$ we have that $(g \circ \psi)^{*}T_{X} \cong g^{*}T_{Y}$ and $N_{g/Y} \cong N_{\psi \circ g/X}$.
\end{lemma}

\begin{proof}
There is an exact sequence $0 \to T_Y \to \psi^*T_X \to T \to 0$, where $T$ is a torsion sheaf.  Taking the top exterior power, we get an injection $\omega_{Y}^{\vee} \to \psi^{*}\omega_{X}^{\vee}$.  We let $R$ denote the ramification divisor whose ideal sheaf is defined by the injection $\omega_{Y}^{\vee} \otimes \psi^{*}\omega_{X} \to \mathcal{O}_{Y}$.  Since $T$ is supported on $\Supp(R)$ we see that $\psi$ is smooth on the complement of $\Supp(R)$. 

Suppose that $C_{X}$ is a general curve parametrized by $M$ and that $C_{Y}$ is a curve parametrized by $N$ that maps to $C_{X}$.  We have
\begin{equation*}
-K_{X} \cdot C_{X} + \dim(X) - 3 = \dim(M) = \dim(N) = -K_{Y} \cdot C_{Y} + \dim(Y) - 3
\end{equation*}
showing that $R \cdot C_{Y} = 0$.  Thus $C_{Y}$ is contained in the smooth locus of $\psi$, so that the map $T_{Y}|_{C_{Y}} \to \psi^{*}T_{X}|_{C_{Y}}$ is surjective.  But a surjective morphism of two locally free sheaves on $\mathbb{P}^{1}$ which have the same rank must be an isomorphism.
\end{proof}

\subsection{Equivalence relation}

\begin{definition}
Let $X$ be a smooth projective variety.  Suppose that $M$ is an irreducible component of $\Rat(X)$ and that we have two $M$ factoring covers $f_{1}: Y_{1} \to X$ and $f_{2}: Y_{2} \to X$ with corresponding irreducible components of curves $N_1$ and $N_2$.   For $i=1,2$ let $\mathcal{C}_{i}$ denote the universal family over $N_{i}$.  We say the two factoring covers are \emph{equivalent} if the general fibers of the evaluation maps $\mathcal{C}_{1} \to Y_{1}$ and $\mathcal{C}_{2} \to Y_{2}$ have the same number of irreducible components.
\end{definition}

We have in mind two key examples.  First, if we modify a factoring cover by a birational morphism (as in Example \ref{exam:birationalfactoringcover}) we obtain an equivalent factoring cover.  Our second key example arises from base change:

\begin{example} \label{exam:basechangefactoringcover}
Suppose that $f: Y \to X$ and an irreducible component $N$ of $\Rat(Y)$ yields a factoring cover for the irreducible component $M$ of $\Rat(X)$.  Suppose furthermore there is a morphism $\pi: Y \to Z$ that contracts the curves parametrized by $N$ so that we get a morphism $N \to Z$.  Take any generically finite map $W \to Z$ and let $Y'$ denote the unique irreducible component of $W \times_{Z} Y$ that dominates both $W$ and $Y$ under the projection maps and let $N'$ be the similarly-defined irreducible component of $W \times_{Z} N$.  Then $Y'$ equipped with $N'$ is a factoring cover for $M$ that is equivalent to our original cover.
\end{example}

\begin{theorem} \label{theo:finiteequivalenceclasses}
Let $X$ be a smooth projective variety.  Suppose that $M$ is a component of $\Rat(X)$. There is a finite set of equivalence classes of factoring covers for $M$.
\end{theorem}

\begin{proof}
We will show that if $f: Y \to X$ and the irreducible component $N \subset \Rat(Y)$ is a factoring cover of $M$ then the number of irreducible components of a general fiber of the evaluation map for $N$ is bounded above by the number of irreducible components of a general fiber of the evaluation map for $M$.  The desired finiteness follows immediately.

Fix a general point $y \in Y$ and let $x = f(y)$.  Let $N_{y}$ denote the sublocus of $N$ parametrizing curves through $y$, and similarly for $M_{x}$ and $x$.  We let $\Gamma_{y}$ denote the set of irreducible components of $N_{y}$ and $\Gamma_{x}$ denote the set of irreducible components of $M_{x}$.  Then the pushforward map induces a function $q: \Gamma_{y} \to \Gamma_{x}$.  An incidence correspondence argument shows that (since $y$ is general) for any fixed irreducible component of $N_{y}$ the general curve parametrized by this component is smooth at $y$.

We claim that this function $q$ is injective.  Since $\dim(Y) = \dim(X)$ and $\dim(N) = \dim(M)$ we see that every irreducible component of $N_{y}$ and every irreducible component of $M_{x}$ has the same dimension $\dim(M) + 1 - \dim(X)$.  In particular, if two components $T_{1},T_{2} \in \Gamma_{y}$ map to the same irreducible component of $\Gamma_{x}$, then for a general rational curve $C_{1}$ parametrized by $T_{1}$ there is a rational curve $C_{2}$ parametrized by $T_{2}$ such that $f(C_{1}) = f(C_{2})$.  This implies that $C_{1} = C_{2}$.  Indeed, as discussed above by generality the curves $C_{1}$, $C_{2}$ are smooth (and hence unibranch) at $y$.  Since $y$ is general we know that $f$ is \'etale on an open neighborhood of $y$.  We conclude that $C_{1}$ and $C_{2}$ must be the same on a formal local neighborhood of $y$, hence everywhere.
Since $C_{1} = C_{2}$ we see that $T_{1} = T_{2}$ and the map $q$ is injective.
\end{proof}

\begin{example} 
Let $X$ be a smooth projective variety and $M$ be an irreducible component of $\Rat(X)$.  Suppose that there is a Galois cover $f: Y \to X$ and an irreducible component $N$ of $\Rat(Y)$ that yield a factoring cover for $M$.  As demonstrated in the proof of Lemma \ref{lemm:tangentbundlefactoringcover}, this implies that the general curve parametrized by $N$ is disjoint from the ramification divisor for $Y$.  Since $f$ is a Galois cover, we deduce that a general curve $C$ parametrized by $M$ is disjoint from the branch divisor for $f$.  In particular, the preimage of $C$ will be a disjoint union of $\deg(f)$ rational curves.

Fix a general point $x \in X$ and let $y_{1},\ldots,y_{d}$ denote the preimages of $x$.  Let $N_{1},\ldots,N_{r}$ be all the irreducible components of $\Rat(Y)$ which map dominantly onto $M$ under the pushforward map.  We let $m(N_{i},y_{j})$ denote the number of components of the sublocus of $N_{i}$ parametrizing rational curves through $y_{j}$, and we define $m(M,x)$ similarly.  Then for any $j = 1,\ldots,d$ we have
\begin{equation*}
\sum_{i} m(N_{i},y_{j})  = m(M,x).
\end{equation*}
Indeed, a general rational curve parametrized by $M$ through $x$ is mapped to by some rational curve through $y_{j}$, and as argued in the proof of Theorem \ref{theo:finiteequivalenceclasses} this implies a local bijection on components of the space of rational curves.

In particular, this shows that $N$ is equivalent to $M$ if and only if it is the unique irreducible component of $\Rat(Y)$ mapping to $M$.
\end{example}

\subsection{Factoring covers and Fano varieties}

Suppose that $X$ is a Fano variety and that $f: Y \to X$ is a factoring cover for an irreducible component of $\Rat(X)$.  Since the ramification divisor for $f$ has vanishing intersection against a dominant family of rational curves on $Y$, it lies on the boundary of the pseudo-effective cone.  \cite{LST18} proves a ``finiteness'' statement for the set of all generically finite covers of $X$ which have this property.  The following result translates this ``finiteness'' to the setting of equivalence classes of factoring covers for components of $\Rat(X)$.

\begin{theorem} \label{theo:finitecovers}
Let $X$ be a smooth projective Fano variety.  There is a finite collection of generically finite morphisms $\psi_{i}: Y_{i} \to X$ such that if an irreducible component $M$ of $\Rat(X)$ admits a factoring cover $\psi: Y \to X$ with an irreducible component $N$ of $\Rat(Y)$ then there will be some index $i$ and a component $N_{i}$ of $\Rat(Y_{i})$ such that $\psi_{i}: Y_{i} \to X$ is a factoring cover for $M$ which is equivalent to $N$.
\end{theorem}

\begin{proof}
Suppose that $f: Y \to X$ equipped with the component $N$ of $\Rat(Y)$ is a factoring cover for an irreducible component $M$ of $\Rat(X)$.  As discussed in the proof of Lemma \ref{lemm:tangentbundlefactoringcover}, this implies that the ramification divisor $R$ for $f$ has vanishing intersection against a general rational curve parametrized by $N$.  In particular, $R = K_{Y} - f^{*}K_{X}$ is on the boundary of the pseudo-effective cone.  Since $-f^{*}K_{X}$ is a big and nef divisor, \cite{BCHM} constructs a canonical model $\pi: Y \dashrightarrow T$ associated to the divisor $K_{Y} -f^{*}K_{X}$.  Note that $\pi$ contracts the general curve parametrized by $N$.  By replacing $Y$ be a birational model (and $N$ by the strict transform family), we may assume that $\pi$ is a morphism.

Let $\{ f_{i}: Y_{i} \to X \}$ denote the finite set of morphisms constructed by \cite[Theorem 1.5]{LST18}.  For each $i$, let $\pi_{i}: Y_{i} \to Z_{i}$ denote the canonical map associated to $K_{Y_{i}} -f_{i}^{*}K_{X}$.  \cite[Theorem 1.5]{LST18} shows that for any factoring cover $f: Y \to X$ there exists a variety $Z$ admitting a generically finite morphism $Z \to T$ such that if $Y_{Z}$ denotes a suitably chosen birational model of the main component of the base change $Y \times_{T} Z$ then the induced morphism $f_{Z}: Y_{Z} \to X$ factors through one of the $f_{i}$.  In fact, the proof shows a little more: since $f_{Z}$ is dominant, a general fiber of $\pi_{Z}: Y_{Z} \to Z$ will map birationally to a general fiber of the map $\pi_{i}$.

Since the curves parametrized by $N$ are in the fibers of $\pi$, by taking a base change we obtain a family of rational curves $N_{Z}$ contained in the fibers of $\pi_{Z}$.  As explained in Example \ref{exam:basechangefactoringcover}, the factoring cover $Y,N$ is equivalent to the factoring cover $Y_{Z},N_{Z}$.  Define $N_{i}$ by taking the images of these curves parametrized by $N_{Z}$ under the rational map to $Y_{i}$.  Then it is clear that $Y_{i},N_{i}$ is a factoring cover of $X,M$.  Furthermore, since the the fibers of $\pi_{Z}$ map birationally to the fibers of $\pi_{i}$ the factoring cover $Y_{i},N_{i}$ is equivalent to the cover $Y_{Z},N_{Z}$ and hence also to the original cover $Y,N$.
\end{proof}

We will also need to know a little more about the nef cone of curves for these covers.

\begin{theorem} \label{theo:relativepolyhedrality}
Let $X$ be a smooth projective Fano variety.  Suppose that $\psi: Y \to X$ is a generically finite morphism such that the divisor $R := K_{Y} - \psi^{*}K_{X}$ is on the boundary of the pseudo-effective cone.  Then the face $F$ of $\Nef_{1}(Y)$ consisting of  curve classes $\alpha$ which satisfy $R \cdot \alpha = 0$ is a polyhedral subcone.
\end{theorem}

\begin{proof}
For any $\epsilon > 0$ define $R_{\epsilon} := K_{Y} - (1-\epsilon) \psi^{*}K_{X}$.  It suffices to show that for some sufficiently small $\epsilon$ there is a finite set of nef curve classes $\{ \alpha_{i} \}_{i=1}^{q}$ such that
\begin{equation*}
\Nef_{1}(X) = \Nef_{1}(X)_{R_{\epsilon} \geq 0} + \sum \mathbb{R}_{\geq 0} \alpha_{i}
\end{equation*}
where $\Nef_{1}(X)_{R_{\epsilon} \geq 0}$ denotes the intersection of $\Nef_{1}(X)$ with the halfspace of classes with non-negative intersection against $R_{\epsilon}$.  This is a consequence of Batyrev's conjecture describing a Cone Theorem for the nef cone of curves when applied to the klt pair $(Y,-\psi^{*}K_{X})$; see \cite[Lemma 6.1]{Araujo10} or \cite[Section 3.2]{LST18}.
\end{proof}

Finally we will need a result of \cite{Neumann10}.

\begin{theorem}[\cite{Neumann10} Proposition 3.3.5]
Let $X$ be a smooth projective variety.  Suppose that $\mathcal{C} \subset \Nef_{1}(X)$ is a polyhedral subcone.  Then $\mathcal{C}$ admits a decomposition into a finite union of locally closed convex cones $\cup_{j} \mathcal{C}_{j}$ such that for any fixed $j$ there is a positive integer $r_{j}$ and a linear function $\widehat{\ell}_{j}: \mathcal{C}_{j} \to \mathbb{Q}^{r}$ such that the Harder-Narasimhan filtration
\begin{equation*}
0 = \mathcal{F}_{0} \subset \mathcal{F}_{1} \subset \ldots \subset \mathcal{F}_{r_{j}} = T_{X}.
\end{equation*}
of $T_{X}$ with respect to any class $\alpha \in \mathcal{C}_{j}$ satisfies that $\mu_{\alpha}(\mathcal{F}_{k}/\mathcal{F}_{k-1})$ is the $k$th entry of $\widehat{\ell}_{j}(\alpha)$.
\end{theorem}

We are now able to prove Theorem \ref{theo:maintheorem2}.

\begin{proof}[Proof of Theorem \ref{theo:maintheorem2}:]
Theorem \ref{theo:finitecovers} shows that there is a finite set of generically finite maps $\psi_{i}: Y_{i} \to X$ such that for every component $M$ of $\Rat(X)$ there will be some index $i$ such that $\psi_{i}$ defines a factoring cover for $M$ that is equivalent to the Stein factorization of the evaluation map for the universal family over $M$.  In particular the corresponding family is an element of $\Rat_{conn}(Y_{i})$.  Furthermore, by Theorem \ref{theo:relativepolyhedrality} there will be a rational polyhedral subcone $\mathcal{C}_{i}$ of $\Nef_{1}(Y_{i})$ which contains the classes of all rational curves on $Y_{i}$ obtained in this way.  Applying \cite[Proposition 3.3.5]{Neumann10} to each variety $Y_{i}$ equipped with the cone $\mathcal{C}_{i}$ we obtain a finite set of linear functions $\widetilde{\ell}_{ij}$ on rational polyhedral subcones of $\mathcal{C}_{i}$.  We obtain the desired functions $\ell_{ij}$ by dividing the $k$th entry of $\widetilde{\ell}_{ij}$ by the rank of the corresponding quotient $\mathcal{F}_{k}/\mathcal{F}_{k-1}$ and repeating this number $\rk(\mathcal{F}_{k}/\mathcal{F}_{k-1})$ times.
\end{proof}

\section{Restricted tangent bundles for Fano varieties}

In this section we return to our original motivating question -- what are the possible restricted tangent bundles for free rational curves on Fano varieties?

\subsection{Slope panels}

Our first collection of results addresses the slope panel of an arbitrary torsion-free sheaf.

\begin{proposition} \label{prop:semistablegm}
Let $X$ be a smooth projective variety and let $\mathcal{E}$ be a torsion-free sheaf on $X$.  Suppose $C$ is a free rational curve that is general in an irreducible component of $\Rat_{conn}(X)$ such that $\mathcal{E}$ is semistable with respect to $[C]$. If we write $\mathcal{E}|_{C} \cong \oplus \mathcal{O}(a_{i})$ then every summand satisfies $|\mu_{C}(\mathcal{E}) - a_{i}| < \frac{\rk(\mathcal{E})}{2}$.
\end{proposition}

\begin{proof}
\cite[Proposition 3.1]{PatelRiedlTseng} shows that $\mathcal{E}|_{C}$ is sequential, and the statement follows immediately.
\end{proof}

\begin{corollary}
Let $X$ be a smooth projective variety and let $\mathcal{E}$ be a torsion-free sheaf on $X$.  Suppose $C$ is a free rational curve that is general in an irreducible component of $\Rat_{conn}(X)$ such that $T_X$ is semistable with respect to $[C]$. Then the minimal slope ratio of $C$ is bounded below by $1-\frac{\dim(X)^{2}}{2(-K_{X} \cdot C)}$.
\end{corollary}

Proposition \ref{prop:semistablegm} can be extended to arbitrary torsion free sheaves using Harder-Narasimhan filtrations.

\begin{corollary} \label{coro:hnfversion}
Let $X$ be a smooth projective variety and let $\mathcal{E}$ be a torsion-free sheaf on $X$ of rank $r$.  Suppose $C$ is a general free curve parametrized by an irreducible component of $\Rat_{conn}(X)$.  Write
\begin{equation*}
0 = \mathcal{F}_{0} \subset \mathcal{F}_{1} \subset \ldots \subset \mathcal{F}_{s} = \mathcal{E}
\end{equation*}
for the Harder-Narasimhan filtration of $\mathcal{E}$.  
Let $\vec{v}$ be the non-increasing $r$-tuple which for every $j \geq 1$ contains $\rk(\mathcal{F}_{j}/\mathcal{F}_{j-1})$ repetitions of the entry $\mu_{C}(\mathcal{F}_{j}/\mathcal{F}_{j-1})$.
If we write $\mathcal{E}|_{C} \cong \oplus \mathcal{O}(a_{i})$ where the $a_{i}$ are in non-increasing order, then 
\begin{equation*}
\left| \vec{v} - \left( a_{1}, \ldots, a_{r} \right) \right|_{sup} < \frac{\sup_{i} \{\rk(\mathcal{F}_{i}/\mathcal{F}_{i-1})\}}{2}.
\end{equation*}
\end{corollary}

\begin{proof}
We prove this statement by induction on the length $s$ of the Harder-Narasimhan filtration.  The base case $s=1$ is Proposition \ref{prop:semistablegm}.

For the induction step, note that since $C$ is general every term $\mathcal{F}_{i}$ is locally free along $C$.  Thus we have an exact sequence
\begin{equation*}
0 \to \mathcal{F}_{s-1}|_{C} \to \mathcal{E}|_{C} \to \mathcal{F}_{s}/\mathcal{F}_{s-1}|_{C} \to 0.
\end{equation*}
We next compare the $\rk(\mathcal{F}_{s-1})$ largest summands in $\mathcal{E}|_{C}$ (written in decreasing order) with the summands of $\mathcal{F}_{s-1}|_{C}$ (in decreasing order).  Any given entry of the former list is no smaller than the corresponding entry on the latter.  The summand for $\mathcal{E}|_{C}$ may be larger, but it will never be larger than the largest summand of $\mathcal{F}_{s}/\mathcal{F}_{s-1}|_{C}$.  This in turn is no larger than
\begin{equation*}
\mu_C(\mathcal{F}_{s}/\mathcal{F}_{s-1}) + \frac{\{\rk(\mathcal{F}_{s}/\mathcal{F}_{s-1})\}}{2}.
\end{equation*}
By induction all the summands of $\mathcal{F}_{s-1}|_{C}$ lie within $\sup_{i \leq s-1} \{\rk(\mathcal{F}_{i}/\mathcal{F}_{i-1})\}/2$ of the expected value.  The argument above shows that the largest $\rk(\mathcal{F}_{s-1})$ direct summands of $\mathcal{E}|_{C}$ can only increase but if they do they will still lie within the desired range.

Similarly, any entry in the list of the $\rk(\mathcal{F}_{s}/\mathcal{F}_{s-1})$ lowest summands of $\mathcal{E}|_{C}$ (written in decreasing order) cannot be larger than the corresponding entry in the list of summands of $\mathcal{F}_{s}/\mathcal{F}_{s-1}|_{C}$ (in decreasing order).  the summand of $\mathcal{E}|_{C}$ can be smaller, but in this case it cannot be smaller than
\begin{equation*}
\mu_C(\mathcal{F}_{s-1}/\mathcal{F}_{s-2}) -  \frac{\sup_{i \leq s-1} \{\rk(\mathcal{F}_{i}/\mathcal{F}_{i-1})\}}{2}.
\end{equation*}
This implies the desired statement for these entries of $\mathcal{E}|_{C}$.
\end{proof}

\subsection{Restricted tangent bundles}

We next turn to the restricted tangent bundle.  If we divide the equation in Corollary \ref{coro:hnfversion} by the slope we obtain:

\begin{corollary} \label{coro:approxforsp}
Let $X$ be a smooth projective variety and fix $\epsilon > 0$.  Let $C$ be a generic curve parametrized by an irreducible component of $\Rat_{conn}(X)$ with anticanonical degree $> \frac{\dim(X)^{2}}{2\epsilon}$.  Then
\begin{equation*}
|\ESP(C) - \SP(C)|_{\sup} < \epsilon.
\end{equation*}
In particular, if $-K_{X}$ is big then there for any $\epsilon$ there are only finitely many families of curves for which this inequality fails to hold.
\end{corollary}

We can translate this into a geometric statement in the following way.

\begin{lemma} \label{lemm:positivityofhn}
Let $X$ be a smooth projective variety.  Let $M$ be an irreducible component of $\Rat_{conn}(X)$ and let $C$ be a general curve parametrized by $M$.  Let
\begin{equation*}
0 = \mathcal{F}_{0} \subset \mathcal{F}_{1} \subset \ldots \subset \mathcal{F}_{r} = T_{X}.
\end{equation*}
denote the Harder-Narasimhan filtration of $T_{X}$ with respect to $C$.  Then for every index $i$ with $1 \leq i \leq r$ we have that $\mu_{C}(\mathcal{F}_{i}/\mathcal{F}_{i-1}) \geq 0$.
\end{lemma}

\begin{proof}
Suppose for a contradiction that $\ESP(C)$ had a negative entry $c$ and fix $\epsilon < |c|$.  Note that if we glue $r$ curves parametrized by $M$ and take a general smoothing of the result to obtain a free curve $C'$ we have $\ESP(C') = \ESP(C)$.  By choosing $r$ sufficiently large, Corollary \ref{coro:approxforsp} guarantees that we can ensure
\begin{equation*}
|\ESP(C) - \SP(C')|_{\sup} < \epsilon.
\end{equation*}
But $\SP(C')$ is non-negative since $C'$ is free.  This yields the desired contradiction.
\end{proof}

\begin{corollary} \label{coro:contractionmaps}
Let $X$ be a smooth projective variety.  Let $M$ be an irreducible component of $\Rat_{conn}(X)$ and let $C$ be a general curve parametrized by $M$.  Let
\begin{equation*}
0 = \mathcal{F}_{0} \subset \mathcal{F}_{1} \subset \ldots \subset \mathcal{F}_{r} = T_{X}.
\end{equation*}
denote the Harder-Narasimhan filtration of $T_{X}$ with respect to $C$.  Then for every index $i$ with $0 < i < r$ the sheaf $\mathcal{F}_{i}$ is a foliation induced by a rational map $\phi: X \dashrightarrow Z$ such that the closures of the fibers of $\phi$ are rationally connected.
\end{corollary}

\begin{proof}
\cite[Theorem 1.1]{CP19} shows that the existence of desired foliation is a consequence of the positivity of the slopes of the quotients in the Harder-Narasimhan filtration as in Lemma \ref{lemm:positivityofhn}.
\end{proof}

Our next goal is to prove Theorem \ref{theo:maintheorem1} by removing the connected fibers condition in Corollary \ref{coro:approxforsp}.  The argument relies on our work with factoring covers, and in particular on Lemma \ref{lemm:tangentbundlefactoringcover} which shows that restricted tangent bundles interact well with factoring covers.  Before proving Theorem \ref{theo:maintheorem1} we need one last definition:

\begin{definition} \label{defi:rgluedcomponents}
Let $X$ be a smooth projective variety.  Fix a positive integer $k$.  For each irreducible component $M$ of $\Rat(X)$, let $\psi: Y \to X$ denote a resolution of the Stein factorization of a normalization of the evaluation map for the universal family over $M$ and let $N$ denote the corresponding irreducible component of $\Rat(Y)$.  By construction the evaluation map for the universal family over $N$ has connected fibers.  Suppose we fix a general point $p$ in $Y$ and take $k$ general curves parametrized by $N$ through $p$.  Since all the curves are free, the resulting stable map represents a smooth point of $\overline{M}_{0,0}(Y)$ and is thus contained in a unique irreducible component.  By applying Lemma \ref{lem-glueNonACovers} repeatedly, we see that this component is independent from the choice of point $p$ and of the general curves through $p$ and thus we may denote it by $N^{(k)}$.  The pushforward map takes $N^{(k)}$ birationally onto an irreducible component of $\Rat(X)$.   We let $\Rat^{(k)}(X)$ denote the union of the components of $\Rat(X)$ obtained in this way as we vary $M$.
\end{definition}

\begin{proof}[Proof of Theorem \ref{theo:maintheorem1}:]
(1) follows immediately from applying Corollary \ref{coro:approxforsp} to a Stein factorization of the evaluation map for the universal family over irreducible components of $\Rat(X)$.  By Lemma \ref{lemm:tangentbundlefactoringcover} the restricted tangent bundle does not change upon this operation.  (2) is a consequence of Theorem \ref{theo:balanceduptorank6} applied to a Stein factorization of the evaluation map for the universal family over irreducible components of $\Rat(X)$.  Again we apply Lemma \ref{lemm:tangentbundlefactoringcover} to deduce that the restricted tangent bundles on $X$ and on the Stein factorization are the same.
\end{proof}

Theorem \ref{theo:maintheorem1} shows that in low dimension for ``sufficiently divisible'' families of curves the slope panel is the same as the expected slope panel.  It is natural to wonder whether ``sufficiently positive'' families of curves will have a slope panel which is ``as close as possible'' to the expected slope panel.  
The following example shows this need not be the case.

\begin{example}
\label{ex-Grassmannian}
Let $\gG(1,3)$ be the Grassmannian of lines in $\PP^3$. Then the tangent bundle $T_{\gG(1,3)}$ is given by $S^{\vee} \otimes Q$, where $S$ and $Q$ are the universal sub and quotient bundles on $\gG(1,3)$. By \cite{Mandal} the restrictions of $S$ and $Q$ to a general rational curve of given degree on $\gG(1,3)$ will be as balanced as possible. Thus, for curves on $\gG(1,3)$ of odd degree, we have that $S^{\vee} \otimes Q$ will have the form $(\OO(a) \oplus \OO(a+1)) \otimes (\OO(a) \oplus \OO(a+1)) = \OO(2a) \oplus \OO(2a+1)^2 \oplus \OO(2a+2)$.  In particular, the restricted tangent bundle will not be $\OO(2a+1)^{4}$ even though there is no arithmetic obstruction to such a decomposition.
\end{example}

\section{Minimal slope ratio} \label{sect:gmc}

Recall that the minimal slope ratio of a free rational curve is the minimal entry in the slope panel.  \cite{Peyre17} proposes a variant of Manin's Conjecture for rational points that depends on an analogue of the minimal slope ratio.  In this section we analyze Peyre's proposal in the setting of Geometric Manin's Conjecture for rational curves.

\subsection{Geometric Manin's Conjecture}

Let $X$ be a smooth projective Fano variety.  Geometric Manin's Conjecture predicts the asymptotic growth rate of the number of components of $\Rat(X)$ as the anticanonical degree increases.  To get the ``expected'' counting function, we must first discount certain families of curves.  One option is given by the following definition.

\begin{definition} \label{defi:pathologicalcomponents}
Let $X$ be a smooth projective Fano variety.  We call an irreducible component of $\Rat(X)$ \emph{pathological} if it admits a factoring cover $f: Y \to X$ such that $f$ is not a birational map.
\end{definition}

\begin{remark}
Any irreducible component of $\Rat(X)$ that fails to be a component of $\Rat_{conn}(X)$ will be pathological, but there can be other pathological components.  For example, if we have a morphism $\pi: X \to Z$ with $\dim Z \geq 1$ then any family of free curves contracted by $\pi$ will be pathological.
\end{remark}

\begin{remark}
The notion of a pathological component is very similar to the notion of an ``accumulating component'' as defined by \cite[Definition 7.1]{BLRT20} which depends upon the exceptional set as constructed in \cite{LST18}.  The two definitions are not precisely the same, but \cite[Theorem 7.7]{BLRT20} shows that for Fano threefolds the distinction between the two is negligible compared to the asymptotic growth of the counting function.
\end{remark}

The counting function in Geometric Manin's Conjecture captures the number of non-pathological components of $\Rat(X)$ of a given degree.  The following version emphasizes the relationship with Manin's Conjecture over a number field.

\begin{definition}
Let $X$ be a smooth Fano variety and let $r(X,-K_{X})$ denote the minimal positive integer of the form $K_{X} \cdot \alpha$ for some $\alpha \in N_{1}(X)_{\mathbb{Z}}$.  We will let $\mathcal{P}_{i}$ denote the set of non-pathological components of $\Rat(X)$ which parametrize curves of anticanonical degree $i \cdot r(X,-K_{X})$

Fix a positive constant $q>1$.  We define
\begin{equation*}
N(X, -K_{X}, q, d) = \sum_{i=1}^{d} \sum_{W \in \mathcal{P}_{i}} q^{\dim W}.
\end{equation*}
\end{definition}

Consider the function $\xi: N_{1}(X)_{\mathbb{Z}} \to \mathbb{Z}_{\geq 0}$ which assigns to any curve class $\alpha$ the number of non-pathological components of $\Rat(X)$ parametrizing curves of class $\alpha$.  Then the counting function is entirely determined by the function $\xi$:
\begin{equation*}
N(X, -K_{X}, q, d) = \sum_{\substack{\alpha \in \Nef_{1}(X)_{\mathbb{Z}} \\ -K_{X} \cdot \alpha \leq d r(X,-K_{X})}} \xi(\alpha) q^{-K_{X} \cdot \alpha}.
\end{equation*}
With this notation, we can state a version of Geometric Manin's Conjecture:

\begin{conjecture} \label{conj:gmc}
Let $X$ be a smooth projective Fano variety.  Then:
\begin{enumerate}
\item There is some constant $M$ such that $\xi(\alpha) \leq M$ for every $\alpha \in N_{1}(X)_{\mathbb{Z}}$.
\item for some translate $\mathcal{Q} = \beta + \Nef_{1}(X)$ of the nef cone, every numerical class in $\mathcal{Q}$ satisfies $\xi(\alpha) = |\Br(X)|$.
\end{enumerate}
\end{conjecture}

Using standard lattice summation techniques, Conjecture \ref{conj:gmc} would imply that the counting function satisfies an asymptotic formula
\begin{equation*}
N(X,-K_{X},q,d) \sim_{d \to \infty} C q^{dr(X, -K_X)} d^{\rho(X)-1}
\end{equation*}
for some constant $C$.

\subsection{Peyre's formulation of Manin's Conjecture}

\begin{definition}
Let $X$ be a smooth projective variety.  Fix $\epsilon > 0$.  We say that a free rational curve $C$ is $\epsilon$-liberated if its minimal slope ratio is $> \epsilon$.
\end{definition}

Suppose that $\epsilon: [1,\infty) \to (0,1)$ is a continuous decreasing function such that $\lim_{d \to \infty} \epsilon(d) = 0$.   By analogy with \cite[D\'efinitions 6.11]{Peyre17}, we define
\begin{equation*}
N^{\ell > \epsilon}(X,-K_{X},q,d) = \sum_{i=1}^{d} \sum_{W \in \widetilde{\mathcal{P}}_{i}} q^{\dim W}.
\end{equation*}
where $\widetilde{\mathcal{P}}_{i}$ denotes the set of $\epsilon(d)$-liberated components of $\Rat(X)$ which parametrize curves of anticanonical degree $i \cdot r(X,-K_{X})$.

\cite{Sawin20} showed that when working over a number field Peyre's definition can fail to remove ``enough'' rational points to obtain the expected growth rate in Manin's Conjecture.  Analogously, in the setting of Geometric Manin's Conjecture there can be pathological components of $\Rat(X)$ which are $\epsilon$-liberated for all sufficiently small $\epsilon$.  Furthermore, such components can make a non-trivial contribution to the asymptotics of the counting function so that $N^{\ell > \epsilon}(X,-K_{X},q,d)$ fails to have the expected growth rate.

\begin{example}[\cite{Sawin20}]
Let $X$ denote the Fano variety $\mathrm{Hilb}^{2}(\mathbb{P}^{n})$.  If we let $\phi: W \to \mathbb{P}^{n} \times \mathbb{P}^{n}$ denote the blow-up of the diagonal, then the quotient of $W$ by the involution that switches the two factors induces a finite map $f: W \to X$ which is ramified along the $\phi$-exceptional divisor $E$.

Suppose that $N$ is an irreducible component of $\Rat(W)$ such that the general curve $C'$ parametrized by $N$ does not intersect $E$.  (We can obtain many such families by taking the strict transforms of irreducible components of $\Rat(\mathbb{P}^{n} \times \mathbb{P}^{n})$.)  If $M$ denotes the irreducible component of $\Rat(X)$ parametrizing the images $C$ of these curves in $X$ then $f$ is a factoring cover for $M$.  In particular, $M$ is a pathological component of $\Rat(X)$ and should not be counted in Geometric Manin's Conjecture.

On the other hand, since a general curve $C'$ parametrized by $N$ avoids $E$ the restriction of the map $T_{W} \to f^{*}T_{X}$ to $C'$ is an isomorphism.  By the same logic, the restriction of the map $T_{W} \to \phi^{*}T_{\mathbb{P}^{n} \times \mathbb{P}^{n}}$ to $C'$ yields an isomorphism.  By choosing different components of $\Rat(\mathbb{P}^{n} \times \mathbb{P}^{n})$, we see that for any fixed $\epsilon$ there will be many components of $\Rat(X)$ coming from $W$ that are pathological but are $\epsilon$-liberated.  Since $\rho(\mathbb{P}^{n} \times \mathbb{P}^{n}) = 2 = \rho(X)$, these curves will affect the leading constant in the asymptotic formula.

(If we are willing to work with weak Fano varieties in place of Fano varieties, then by doing a similar construction for $X = \Hilb^{2}(\PP^{m} \times \PP^{m})$ we can even ensure that the exponents in the asymptotic formula for $W$ are larger than the exponents for $X$.)
\end{example}

We next show that Peyre's proposal never removes ``too much'': although Peyre's suggestion will discount some non-pathological families, if we assume Conjecture \ref{conj:gmc} such families do not impact the asymptotic formula.  We will need the following facts:

\begin{observation} \label{obse:slopesandcones}
\begin{enumerate}
\item  By taking closures of the cones constructed in \cite[Proposition 3.3.5]{Neumann10}, the cone $\Nef_{1}(X)$ can be decomposed into a finite set of closed full-dimensional rational polyhedral cones $\{ \mathcal{C}_{i} \}_{i=1}^{r}$ such that for every $i$ the Harder-Narasimhan filtration of $T_{X}$ with respect to classes $\alpha$ in the interior of $\mathcal{C}_{i}$ is constant.
\item Conjecture \ref{conj:gmc}.(2) implies that the closure of the cone generated by the classes of components of $\Rat_{conn}(X)$  is all of $\Nef_{1}(X)$.  By Lemma \ref{lemm:positivityofhn} this implies that for any $\mathcal{C}_{i}$ and any successive quotient $\mathcal{F}_{k}/\mathcal{F}_{k-1}$ occurring in the corresponding Harder-Narasimhan filtration the slope of $\mathcal{F}_{k}/\mathcal{F}_{k-1}$ with respect to any numerical class in $\mathcal{C}_{i}$ is non-negative.  Since $\mathcal{C}_{i}$ is full-dimensional it is not possible for these slopes to be identically $0$ on $\mathcal{C}_{i}$, and thus each will be positive on the interior of $\mathcal{C}_{i}$.
\end{enumerate}
\end{observation}

\begin{theorem} \label{theo:peyreexceptionalset}
Assume Conjecture \ref{conj:gmc}.  Let $X$ be a smooth projective Fano variety and let $\epsilon: [1,\infty) \to (0,1)$ be a continuous decreasing function such that $\lim_{d \to \infty} \epsilon(d) = 0$.  For any $\delta > 0$, there is some $d_{0}$ sufficiently large such that
\begin{equation*}
N^{\ell > \epsilon}(X,-K_{X},q,d) > (1-\delta) N(X,-K_{X},q,d)
\end{equation*}
for every $d \geq d_{0}$.
\end{theorem}

\begin{proof}
We continue to write $\{ \mathcal{C}_{i} \}_{i=1}^{b}$ for the subcones of $\Nef_{1}(X)$ such that the Harder-Narasimhan filtration is constant for classes in the interior of $\mathcal{C}_{i}$.  We let $\ell_{i}$ denote the linear function on $\mathcal{C}_{i}$ that identifies the minimal slope quotient.

Let $\mathcal{Q} \subset \Nef_{1}(X)$ denote the translate of the nef cone identified by Conjecture \ref{conj:gmc}.(2).  Conjecture \ref{conj:gmc}.(1) implies that for $d$ sufficiently large we have
\begin{equation} \label{eq:latticecount}
\sum_{\substack{\alpha \in \mathcal{Q}_{\mathbb{Z}} \\ -K_{X} \cdot \alpha \leq d r(X,-K_{X})}} \xi(\alpha) q^{-K_{X} \cdot \alpha} > (1-\delta)^{1/6} N(X, -K_{X}, q, d).
\end{equation}
Furthermore Conjecture \ref{conj:gmc}.(2) predicts that $\xi(\alpha) = |\Br(X)|$ for every $\alpha \in \mathcal{Q}$.  Together these imply that we can approximate $N(X,-K_{X},q,d)$ using a sum of exponentials over the lattice points of $\mathcal{Q}$.  Our plan is to estimate the percentage of numerical classes in $\mathcal{Q}$ which represent $\epsilon$-liberated rational curves.

For any positive integer $d$ let $\mathcal{N}^{d}$ denote the subset in $\Nef_{1}(X)$ of classes which have anticanonical degree $= r(X,-K_{X}) \cdot d$.  By Observation \ref{obse:slopesandcones}.(2) we see that for $\gamma$ sufficiently small we have
\begin{equation*}
\Vol \left( \left\{ \alpha \in \mathcal{N}^{d} \cap \mathcal{C}_{i} \left| \ell_{i}(\alpha) \geq \frac{\gamma}{2}(-K_{X} \cdot \alpha) \right. \right\} \right) > (1-\delta)^{1/6} \Vol(\mathcal{N}^{d} \cap \mathcal{C}_{i})
\end{equation*}
where the volume is normalized with respect to the lattice of integer curve classes in the codimension $1$ subspace of $N_{1}(X)$ consisting of classes with $-K_{X} \cdot \alpha = 0$.  Appealing to the theory of Ehrhart quasipolynomials, we see that by decreasing the leading constant we may ensure that for any fixed value of $\gamma$ as above there is a constant $d_{0}$ such that
\begin{equation*}
\# \left\{ \alpha \in \mathcal{N}_{\mathbb{Z}}^{d} \cap \mathcal{C}_{i} \left| \ell_{i}(\alpha) \geq \frac{\gamma}{2}(-K_{X} \cdot \alpha) \right. \right\} > (1-\delta)^{2/6} \#(\mathcal{N}_{\mathbb{Z}}^{d} \cap \mathcal{C}_{i})
\end{equation*}
for every $d \geq d_{0}$.  

Since $\mathcal{Q}$ is a translate of the nef cone, for sufficiently large $d$ the number of lattice points in a cross section of $\mathcal{Q} \cap \mathcal{C}_{i} \cap \mathcal{N}_{\mathbb{Z}}^{d}$ is asymptotically the same as the number in $\mathcal{C}_{i} \cap \mathcal{N}_{\mathbb{Z}}^{d}$.  More precisely, after perhaps increasing $d_{0}$ further to account for the translate defining $\mathcal{Q}$ (and increasing the leading constant) we see that for every $d \geq d_{0}$ we have
\begin{equation*}
\# \left\{ \alpha \in \mathcal{Q} \cap \mathcal{N}_{\mathbb{Z}}^{d} \cap \mathcal{C}_{i} \left| \ell_{i}(\alpha) \geq \frac{\gamma}{2}(-K_{X} \cdot \alpha) \right. \right\} > (1-\delta)^{3/6} \#(\mathcal{Q} \cap \mathcal{N}_{\mathbb{Z}}^{d} \cap \mathcal{C}_{i}).
\end{equation*}
Since ``almost all'' of the points in $\mathcal{Q} \cap \mathcal{C}_{i} \cap \mathcal{N}_{\mathbb{Z}}^{d}$ satisfy the inequality $\ell_{i}(\alpha) \geq \frac{\gamma}{2}(-K_{X} \cdot \alpha)$ for $d$ sufficiently large, we see that after possibly increasing $d_{0}$ (and the leading constant) we have for every $d \geq d_{0}$
\begin{equation*}
\sum_{\substack{\alpha \in  \mathcal{Q} \cap \mathcal{N}_{\mathbb{Z}}^{d} \cap \mathcal{C}_{i} \\ \ell_{i}(\alpha) \geq \frac{\gamma}{2}(-K_{X} \cdot \alpha) }} q^{-K_{X} \cdot \alpha} > (1-\delta)^{4/6} \sum_{\substack{\alpha \in  \mathcal{Q} \cap \mathcal{N}_{\mathbb{Z}}^{d} \cap \mathcal{C}_{i} }} q^{-K_{X} \cdot \alpha}.
\end{equation*}
Summing up over anticanonical degrees, we see that for $d_{0}$ sufficiently large we have
\begin{equation}  \label{eq:gammadeltaineq}
\sum_{\substack{\alpha \in  \mathcal{Q}_{\mathbb{Z}} \cap \mathcal{C}_{i} \\ \\   -K_{X} \cdot \alpha \leq d r(X,-K_{X}) \\ \ell_{i}(\alpha) \geq \frac{\gamma}{2}(-K_{X} \cdot \alpha) }} q^{-K_{X} \cdot \alpha} > (1-\delta)^{5/6} \sum_{\substack{\alpha \in  \mathcal{Q}_{\mathbb{Z}} \cap \mathcal{C}_{i} \\  -K_{X} \cdot \alpha \leq d r(X,-K_{X})  }} q^{-K_{X} \cdot \alpha}.
\end{equation}
for every $d \geq d_{0}$.
In particular we can first choose $\gamma$ small enough and then choose $d_{0}$ sufficiently large so that this inequality holds for every $\mathcal{C}_{i}$ simultaneously.

Suppose that $M$ is an irreducible component of $\Rat_{conn}(X)$ and that $C$ is a general free rational curve parametrized by $M$.  Suppose furthermore that $-K_{X} \cdot C > \dim(X)^{2}/4\gamma$.  Then by Theorem \ref{theo:maintheorem1} we see that if the numerical class of $C$ lies in $\mathcal{C}_{i}$ and if
\begin{equation*}
\ell_{i}(C) \geq \frac{\gamma}{2} (-K_{X} \cdot C)
\end{equation*}
then $C$ is $\gamma$-liberated.  In particular, if $\alpha$ is a numerical class counted by the left hand side of Equation \eqref{eq:gammadeltaineq} then any free rational curve of class $\alpha$ that is general in its deformation class will be $\gamma$-liberated.

Altogether, if we further increase $d_{0}$ to ensure that $\epsilon(d) < \gamma$ for every $d \geq d_{0}$, then we see that for $d \geq d_{0}$ we have
\begin{align*}
N^{\ell > \epsilon}(X,-K_{X},q,d) & \geq \sum_{i=1}^{b} \sum_{\substack{\alpha \in  \mathcal{Q}_{\mathbb{Z}} \cap \mathcal{C}_{i} \\   -K_{X} \cdot \alpha \leq d r(X,-K_{X}) \\ \ell_{i}(\alpha) \geq \frac{\gamma}{2}(-K_{X} \cdot \alpha) }}|\Br(X)| q^{-K_{X} \cdot \alpha}  & \\
& > (1-\delta)^{5/6} \sum_{i=1}^{b} \sum_{\substack{\alpha \in  \mathcal{Q}_{\mathbb{Z}} \cap \mathcal{C}_{i}  \\  -K_{X} \cdot \alpha \leq d r(X,-K_{X})}} |\Br(X)| q^{-K_{X} \cdot \alpha} &  \textrm{by Equation \eqref{eq:gammadeltaineq}} \\
&= (1-\delta)^{5/6}  \sum_{\substack{\alpha \in \mathcal{Q}_{\mathbb{Z}} \\  -K_{X} \cdot \alpha \leq d r(X,-K_{X})}} |\Br(X)| q^{-K_{X} \cdot \alpha}  &\\
& > (1-\delta) N(X, -K_{X}, q, d) & \textrm{by Equation \eqref{eq:latticecount}}
\end{align*}
\end{proof}

\bibliographystyle{amsalpha}
\bibliography{restrictedtangentbundles}

\end{document}